\documentclass[10pt]{article}
\usepackage{graphicx}
\usepackage{amsmath,amssymb,amsthm,amsfonts}
\usepackage{color}
\usepackage{amssymb}
\usepackage{cite}
\newtheorem{thm}{Theorem}[section]

\newtheorem{lem}[thm]{Lemma}
\newtheorem{prop}[thm]{Proposition}

\theoremstyle{definition}

\theoremstyle{remark}
\usepackage{appendix}

\numberwithin{equation}{section}

\DeclareMathSymbol{\C}{\mathalpha}{AMSb}{"43}

\newcommand{\eps}{\varepsilon}

\newcommand{\R}{{\mathbb{R}}}

\def\R{{\mathbb R}}
\def\C{{\mathbb C}}

\newcommand{\bsub}{\begin{subequations}}
\newcommand{\esub}{\end{subequations}$\!$}

\begin{document}

\title{Ground States of Attractive Bose Gases in Rotating Anharmonic Traps}
\date{\today}
\author{Yujin Guo\thanks{School of Mathematics and Statistics, and Hubei Key Laboratory of
Mathematical Sciences, Central China Normal University, P.O. Box 71010, Wuhan 430079, P. R. China. Email: \texttt{yguo@ccnu.edu.cn}. Y. J. Guo is partially supported by NSFC under Grants 12225106 and 11931012.
},
\, 	Yan Li\thanks{School of Mathematics and Statistics, Central China Normal University, P.O. Box 71010, Wuhan 430079,
	P. R. China.  Email: \texttt{yanlimath@mails.ccnu.edu.cn}.},
\, and\, Yong Luo\thanks{School of Mathematics and Statistics, and Hubei Key Laboratory of
		Mathematical Sciences, Central China Normal University, P.O. Box 71010, Wuhan 430079,
		P. R. China.  Email: \texttt{yluo@ccnu.edu.cn}.}
	}

\smallbreak \maketitle

\begin{abstract}
This paper is concerned with ground states of attractive Bose gases confined in an anharmonic trap $V(x)=\omega(|x|^2+k|x|^4)$ rotating at the velocity $\Omega>0$, where $\omega>0$ denotes the trapping frequency, and $k>0$ represents the strength of the quartic term. It is known that for any $\Omega>0$, ground states exist in such traps if and only if $0<a<a^*$, where $a^*:=\|Q\|^{2}_{2}$ and $Q>0$ is the unique positive solution of $\Delta Q-Q+Q^{3}=0$ in $\R^2$. By analyzing the refined energies and expansions of ground states, we prove that there exists a constant $C>0$, independent of $0<a<a^*$, such that ground states do not have any vortex in the region $R(a):=\big\{x\in\R^2:\, |x|\leq C(a^*-a)^{-\frac{1-6\beta}{20}}\big\}$ as $a\nearrow a^*$, for the case where  $\omega=\frac{3\Omega^2}{4}$, $k=\frac{1}{6}$, and $\Omega=C_0(a^*-a)^{-\beta}$ varies for some $\beta\in [0,\frac{1}{6})$ and $C_0>0$.
\end{abstract}



\section{Introduction}

The Bose gases trapped in a magnetic potential present the remarkable property
that a macroscopic fraction of the atoms at sufficiently low temperatures are condensed in
the same quantum state. This phenomenon is referred to as a Bose-Einstein condensate (BEC)
and was first observed experimentally in 1995, see \cite{AEM,DAV}. The particularly
interesting experiments of creating BECs are to
cool atomic gases confined in a magnetic trap rotating at the velocity $\Omega$.
Various interesting quantum phenomena
have been observed in the experiments of rotating  Bose gases since the late 1990s, including the critical-mass collapse \cite{BSTH,DGPS,HMDB}, the center-of-mass rotation \cite{AS,F,LCS,WGS}, the effective lower dimensional behavior in strongly elongated traps \cite{AA,AEM,BC,DGPS,F}, and so on. Therefore, mathematical theories of rotating Bose gases have become an important topic of researches over the past two decades, see \cite{BLA,GLY,Lewin,LNR,RN,RN1,MCW,RAV} and the references therein.

When the interactions between Bose gases are repulsive, the quantized vortices and some other complex structures of Bose gases confined in rotating traps were analyzed and simulated extensively in the past few years, see \cite{AA1,Lieb2,RN1,IM-1,SSbook} and the references therein. Because of the distinct  critical-mass collapse, the systems of attractive Bose gases  confined in rotating traps however behave quite different from those of the well-understood repulsive case. For example, the vortices are generally unstable in attractive Bose gases  under rotation (see, e.g., \cite{CC,LCS,WGS}), and meanwhile it is known (cf. \cite{AA,Cooper,F}) that the vortices  form
stable lattice configurations in the repulsive case.

As derived rigorously by a mean-field approximation (cf. \cite{Lewin,LNR,Lieb2,NAM,RN1}), ground states of two-dimensional attractive Bose gases confined in a rotating trap can be described equivalently (see \cite{F,LNR,Lieb3,PLP}) by  the  minimizers of the following complex-valued variational problem:
\begin{equation*}
\tilde{I}(a,\Omega):=\inf_{\{u\in\mathcal{H},\ \|u\|_{2}^2=1\}}\tilde{F}_{a,\Omega}(u),\, \ a>0,\,\ \Omega >0,
\end{equation*}
where the energy functional $\tilde{F}_{a,\Omega}(u)$ is given by
\begin{equation}\label{GP}
\begin{aligned}
    \tilde{F}_{a,\Omega}(u)
    &=\int_{\R^2}\big|\big(\nabla-i\Omega x^{\bot}\big)u\big|^2dx+\int_{\R^2}\big(V(x)-\Omega^2|x|^2\big)|u|^2dx-\frac{a}{2}\int_{\R^2}|u|^4dx,
\end{aligned}
\end{equation}
and the space $\mathcal{H}$ is defined as
\begin{equation*}
    \mathcal{H}:=\Big\{u\in H^{1}(\R^2,\C):\ \int_{\R^2}V(x)|u|^2dx<\infty\Big\}.
\end{equation*}
Here $x^{\bot}=(-x_{2},x_{1})$ holds for $x=(x_{1},x_{2})\in\R^2$, $V(x)\ge 0$ is the trapping potential confining cold atoms and rotating at the velocity $\Omega\geq0$. However, the parameter $a>0$ characterizes the absolute product of the scattering length $\nu$ of the two-body interaction times the number $N$ of particles in the condensates.

When there is no rotation for the trap $V(x)$, i.e., $\Omega=0$, the existence, uniqueness, mass concentration and some other analytical properties of  minimizers for $\tilde{I}(a,\Omega)$ were studied recently in \cite{GLW,GS,Guo1,LPWY,ZJ} and the references therein. Specially, it was shown implicitly in above mentioned works that $\tilde{I}(a,0)$ admits minimizers if and only if $a<a^*$, where $a^*:=\|Q\|_{2}^2$ and $Q(x)=Q(|x|)>0$ is the unique positive
solution of the following nonlinear equation
\begin{equation}\label{Q}
    -\Delta Q+Q-Q^3=0\ \ \hbox{in}\ \ \R^2,\ \ u\in H^{1}(\R^2,\R).
\end{equation}
The rotational case $\Omega>0$ of $\tilde{I}(a,\Omega)$ was studied more recently, see the earlier works \cite{ANS,BC,LNR} and the references therein. However, most existing analysis of  $\tilde{I}(a,\Omega)$ mainly focused on the harmonic trapping potentials of the type
\begin{equation}\label{0:3}
    V(x)=x_{1}^{2}+\Lambda^2 x_{2}^{2},\ \ \Lambda>0,
\end{equation}
where $x=(x_{1},x_{2})\in\R^2$. Under the assumption (\ref{0:3}),  $\tilde{I}(a,\Omega)$ admits a finite critical rotational speed $\Omega ^*:=\min\{1, \Lambda\}>0$, which can be defined as
\begin{equation}
	\Omega ^*:=\sup \Big\{\Omega >0:\, \  V(x)- \Omega ^2 |x|^2  \to\infty \,\ \mbox{as}\,\
	|x|\to\infty \Big\}.  \label{Omega}
\end{equation}
It was proved in \cite[Theorem 1.1]{GLY} that for any $\Omega >\Omega ^*$, there is no minimizer of $\tilde{I}(a,\Omega)$; for any $0<\Omega <\Omega ^*$, $\tilde{I}(a,\Omega)$ admits minimizers if and only if $0<a<a^*:=\|Q\|_{2}^2$.
Furthermore, the  refined analytical properties of minimizers for $\tilde{I}(a,\Omega)$ were studied in \cite{ANS,BC,LNR,GLY,Dinh,G} under the assumption (\ref{0:3}), including the uniqueness and stability, the nonexistence of vortices, and so on. Under the harmonic trap (\ref{0:3}),  we remark from above that
we are unable to take an arbitrarily large rotational velocity $\Omega$. This corresponds physically to the fact that once $\Omega>\Omega ^*$, the centrifugal force of the harmonic trap (\ref{0:3}) overcomes its magnetic trapping force.

Due to the above restriction, Fetter proposed in \cite{F1} the stiffer anharmonic trap of the form
\begin{equation}\label{poten1}
    V(x)=\omega(|x|^2+k|x|^4),\ \ \hbox{where}\ \ \omega>0,\ k>0,
\end{equation}
instead of the harmonic one (\ref{0:3}), where $\omega>0$ denotes the trapping frequency, and $k>0$ represents the strength of the quartic term.
The nice feature of the trapping potential (\ref{poten1}) lies in the fact that the centrifugal force is always compensated by the trapping force, and it thus follows theoretically from  (\ref{Omega}) that one can take an arbitrarily large rotational velocity. Anharmonic trapping potentials of the type \eqref{poten1} were achieved experimentally in \cite{BSS,BSC}, where  the quadratic trapping potential was modified by superimposing a blue detuned laser beam to the magnetic trap holding the atoms, see also \cite{AD,BLA,DAN,FJ,FU,JAC,JKL,KAS,KAV,KIM,LUN,RN} for the theoretical and numerical investigations.


Motivated by above facts, we consider the anharmonic trap (\ref{poten1}) with $\omega=\frac{3\Omega^2}{4}$ and $k=\frac{1}{6}$, so that \eqref{poten1} is of the form
\begin{equation}\label{1.01}
    V(x)=\frac{3\Omega^2}{4}\big(|x|^2+\frac{1}{6}|x|^4\big),\ \ \Omega>0,
\end{equation}
where  $\Omega$ denotes the rotational velocity of the trap.
Under the trapping potential \eqref{1.01}, the energy functional $\tilde{F}_{a,\Omega}(u)$ is then reduced to the form
\begin{equation*}
\begin{aligned}
    \tilde{F}_{a,\Omega}(u)
    &=\int_{\R^2}\big|\big(\nabla-i\Omega x^{\bot}\big)u\big|^2dx+\frac{\Omega^2}{8}\int_{\R^2}\big(|x|^4-2|x|^2\big)|u|^2dx
    -\frac{a}{2}\int_{\R^2}|u|^4dx,
\end{aligned}
\end{equation*}
and the associated critical velocity $\Omega ^*$ defined in (\ref{Omega}) satisfies $\Omega ^*=\infty$. Because of the $L^2-$constraint condition, the analysis of the variational problem $\tilde{I}(a,\Omega)$ is thus equivalent to the following form
\begin{equation*}
    I(a,\Omega):=\inf_{\{u\in\mathcal{H},\ \|u\|_{2}^2=1\}}F_{a,\Omega}(u), \,\  a>0,\,\   \Omega >0,
\end{equation*}
where
\begin{equation*}
\begin{aligned}
    F_{a,\Omega}(u)
    &=\int_{\R^{2}}\big|(\nabla-i\Omega x^{\bot})u\big|^{2}dx+\frac{\Omega^{2}}{8}\int_{\R^{2}}\big(|x|^{2}-1\big)^{2}|u|^{2}dx-\frac{a}{2}\int_{\R^{2}}|u|^{4}dx,\\
\end{aligned}
\end{equation*}
and
\begin{equation*}
    \mathcal{H}:=\Big\{u\in H^{1}(\R^2,\C):\int_{\R^2}|x|^4|u|^{2}dx<\infty\Big\}.
\end{equation*}
Recall from \cite{Lieb}  the following diamagnetic inequality:
\begin{equation}\label{dia}
\big|\big(\nabla-i\Omega x^{\bot}\big)u\big|^2\ge   \big|\nabla|u|\big|^{2},\ \ u\in H^1(\R^2,\C).
\end{equation}
By applying  (\ref{dia}) and the variational methods, we therefore have the following existence and non-existence results of $I(a,\Omega)$.

\vskip 0.05truein
\noindent{\textbf{Theorem A (\cite[Theorem 1.1]{GLY}).}}
Let $Q=Q(|x|)>0$ be the unique positive solution of (\ref{Q}). Then for any $\Omega>0$, we have
\begin{enumerate}
  \item If $0<a<a^*:=\|Q\|_{2}^2$, then there exists at least one minimizer of $I(a,\Omega)$;
  \item If $a\geq a^*$, then there is no minimizer of $I(a,\Omega)$.
\end{enumerate}

By the variational argument, if $u_{a}$ is a complex-valued minimizer of $I(a,\Omega)$, then there exists a Lagrange multiplier $\mu_{a}=\mu(a, u_a)\in\R$ such that $u_{a}$ is a normalized complex-valued solution of the following elliptic equation
\begin{equation}\label{ELE1}
    -\Delta u_{a}+2i\Omega (x^{\bot}\cdot\nabla u_{a})+\Omega^{2}|x|^{2}u_{a}+\frac{\Omega^{2}}{8}(|x|^{2}-1)^{2}u_{a}
    =\mu_{a}u_{a}+a|u_{a}|^{2}u_{a}\ \ \hbox{in}\ \ \R^2.
\end{equation}
The nonexistence of vortices for rotating Bose gases in the repulsive case under the anharmonic trap \eqref{poten1} was studied earlier in \cite{AAL,AJR,BSS,KAS,SBCD} and the references therein, where the jacobian estimates, vortex ball constructions, and some other methods were employed. However, as far as we know, the arguments mentioned above are not applicable to our attractive case, due to the non-convex nonlinearity of the energy $ F_{a,\Omega}(u)$. To the best of our knowledge, there are rarely results on the non-existence of vortices for rotating Bose gases in the attractive case, except the recent works \cite{GLY,G} which were however concerned with the cases where the traps are harmonic. More precisely, if the harmonic trap  $V(x)=|x|^2$ is radially symmetric, then the non-existence of vortices for minimizers was proved in \cite{GLY} by developing the method of inductive symmetry. Further, if the harmonic trap $V(x)=x_{1}^2+\Lambda^2x_{2}^2$, where $0<\Lambda<1$, is non-radially symmetric, then  the non-existence of vortices for minimizers was studied in \cite{G} by analyzing the refined expansions of the minimizers.

By considering the anharmonic trap \eqref{1.01}, instead of the harmonic one, the main purpose of this paper is to establish the following non-existence of vortices for minimizers of $I(a,\Omega)$ as $a\nearrow a^*$.


\begin{thm}\label{th4}
Let $u_{a}$ be a complex-valued minimizer of $I(a,\Omega_{a})$ for $0<a<a^*$, where $\Omega_{a}>0$ satisfies
\begin{equation}\label{0:1}
    \Omega_{a}:=C_{0}(a^*-a)^{-\beta}>0,\ \ 0< a<a^*,\ \beta\in [0,\frac{1}{6})\ \ \hbox{and}\ \ C_{0}>0.
\end{equation}
Then there exists a constant $C>0$, independent of $0<a<a^*$, such that
\begin{equation*}
\begin{aligned}
    &|u_{a}|>0\ \ \hbox{in the large region}\ \ R(a): =\big\{x\in\R^2:\ |x|\leq C(a^*-a)^{-\frac{1-6\beta}{20}}\big\} \ \ \hbox{as}\ \ a\nearrow a^*,
\end{aligned}
\end{equation*}
i.e., $u_{a}$ does not admit any vortex in the region $R(a)$ as $a\nearrow a^*$.
\end{thm}

\noindent As far as we know, the existing analysis on the vortices of minimizers $u_{a}$ for $\tilde I(a,\Omega)$ was mainly devoted to the case where the velocity $\Omega>0$ is fixed, e.g. \cite{G,GLY}. Besides this case, our Theorem \ref{th4} also covers the case where the velocity $\Omega _a>0$ satisfies $\Omega _a\nearrow\infty$. Theorem \ref{th4} shows that the region $R(a)$, where the minimizer $u_{a}$ does not admit any vortex, becomes smaller once the rotational velocity $\Omega_{a}$ increases, $i.e.,$ the parameter $\beta \ge 0$ increases. On the other hand, even though the trap \eqref{1.01} of Theorem \ref{th4} is radially symmetric, one cannot follow the method of inductive symmetry developed in \cite{GLY} to analyze the non-existence of vortices for minimizers in the whole plane.


\subsection{Proof strategy of Theorem \ref{th4}}

The purpose of this subsection is to introduce the general strategy of proving Theorem \ref{th4}. Roughly speaking, the proof of Theorem \ref{th4} makes full use of the refined expansions of minimizers $u_{a}$ for $I(a,\Omega_{a})$ as $a\nearrow a^*$, where $\Omega_{a}>0$ satisfies \eqref{0:1}.

The first step of proving Theorem \ref{th4} is to derive the following  $L^{\infty}-$uniform convergence of minimizers $u_{a}$ after rescaling and translation.

\begin{thm}\label{th1}
Let $u_{a}$ be a complex-valued minimizer of $I(a,\Omega_{a})$ for $ 0< a<a^*$, where $\Omega_{a}>0$ satisfies \eqref{0:1}, i.e., $\Omega_{a}=C_{0}(a^*-a)^{-\beta}$ holds for some $\beta\in [0,\frac{1}{2})$ and $C_{0}>0$. Then for any given sequence $\{a_{k}\}$ satisfying $a_{k}\nearrow a^*$ as $k\rightarrow\infty$, there exists a subsequence, still denoted by $\{a_{k}\}$, of $\{a_{k}\}$ such that $u_{a_{k}}$ satisfies
\begin{equation*}
\begin{split}
    &w_{a_{k}}(x):=\frac{(a^*-a_{k})^{\frac{1+2\beta}{4}}}{\sqrt{C_{0}}\lambda}u_{a_{k}}
    \Big(\frac{(a^*-a_{k})^{\frac{1+2\beta}{4}}}{\sqrt{C_{0}}\lambda}x+x_{a_{k}}\Big)
    e^{-i\Big(\frac{(a^*-a_{k})^{\frac{1+2\beta}{4}}}{\sqrt{C_{0}}\lambda}\Omega_{a_{k}}x\cdot x_{a_{k}}^{\bot}-\theta_{a_{k}}\Big)}\rightarrow \frac{Q(x)}{\sqrt{a^*}}\\
\end{split}
\end{equation*}
uniformly in $L^{\infty}(\R^2,\C)$ as $k\to\infty$, where $\theta_{a_{k}}\in[0,2\pi)$ is a suitable constant phase, and $\lambda>0$ is defined by
\begin{equation}\label{lmd}
    \lambda:=\Big[\frac{5}{4}\int_{\R^2}|x|^{2}Q^{2}dx\Big]^{\frac{1}{4}}>0.
\end{equation}
Here $x_{a_{k}}$ is the  unique maximum point of $|u_{a_{k}}|$ satisfying $\lim\limits_{k\to\infty}|x_{a_{k}}|=|x_0|=1$ and
\begin{equation}\label{M:lmd}
    \lim\limits_{k\to\infty}\frac{|x_{a_{k}}|^2-1}{(a^*-a_{k})^{\frac{1+2\beta}{4}}}=0.
\end{equation}
\end{thm}



We note that most of the existing works (e.g.\cite{ANS,GLP1,GLY,LNR}) were focussed on the limiting behavior of minimizers for $\tilde{I}(a,\Omega)$ with a fixed rotational velocity $\Omega>0$. Theorem \ref{th1} addresses the interesting problem imposed in \cite[Remark 2.2]{LNR}, which is concerned with the limiting behavior of minimizers for $I(a,\Omega_{a})$ as $\Omega_{a} >0$ varies. In this case, we notice from Theorem \ref{th1} that the blow-up rate of minimizers depends on the varying rotational velocity $\Omega _a>0$. On the other hand, even though the $L^{\infty}-$uniform convergence of minimizers $u_a$ for $\tilde{I}(a,\Omega)$ was studied recently in \cite{GLP1,GLY,LNR} and the references therein, there exist some extra difficulties in the proof of Theorem \ref{th1}. Firstly, since the rotational velocity $\Omega_{a}$ satisfies $\Omega_{a}\nearrow\infty$ as $a\nearrow a^*$, one cannot use $\bar{\epsilon}_{a}^{-2}:=\int_{\R^2}\big|\nabla u_{a}\big|^2dx$ as in \cite{GLP1,GLY} to estimate the Lagrange multiplier $\mu_{a}$ of \eqref{ELE1}. To overcome this difficulty, we shall define $\epsilon_{a}^{-2}:=\int_{\R^2}\big|\nabla |u_{a}|\big|^2dx$  in a slightly different way. Secondly, the above definition of $\epsilon_{a}$ however yields that one cannot derive directly the limiting behavior of $u_{a}$ as $a\nearrow a^*$. For this reason, stimulated by \cite{Guo2}, we shall first establish the asymptotic behavior of $|u_{a}|$ as $a\nearrow a^*$. Following it, we then get some refined estimates of $I(a,\Omega_{a})$ as $a\nearrow a^*$, based on which we shall be able to study the refined estimates of $u_{a}$ as $a\nearrow a^*$.


As the second step of proving Theorem \ref{th4}, we follow Theorem  \ref{th1} to derive the refined expansions of minimizers $u_{a_{k}}$ as $k\to\infty$. For convenience, we denote $\Psi_{i}\in C^2(\R^2)\cap L^{\infty}(\R^2)$ to be the unique solution of the following problem
\begin{equation}\label{0.1}
\nabla\Psi_{i}(0)=0,\ \ (-\Delta+1-3Q^2)\Psi_{i}(x)=f_{i}(x)\ \ \hbox{in}\, \ \R^2,\ \ i=1,2,
\end{equation}
where $f_{i}(x)$ satisfies
\begin{equation*}\arraycolsep=1.5pt
f_i(x)=\left\{\begin{array}{lll}
&- \displaystyle\Big[|x|^2+\frac{1}{2}(x\cdot x_{0})^{2}\Big]Q(x),  \  \ &\mbox{if}\ \ i=1;\\[3mm]
&-\frac{1}{2}\big(|x|^2+\tilde{C}\big)(x\cdot x_{0})Q(x), \,\ &\mbox{if}\ \ i=2.
\end{array}\right.\end{equation*}
Here the point $x_0\in\R^2$ satisfying $\lim\limits_{k\to\infty}x_{a_{k}}=x_0$ is as in Theorem \ref{th1}, and the constant $\tilde{C}=-\frac{8\lambda^4}{5a^*}<0$ holds for $\lambda >0$ defined by \eqref{lmd}. We also denote $\Phi_{I}\in C^2(\R^2)\cap L^{\infty}(\R^2)$ to be the unique solution of
\begin{equation}\label{0.3}
(-\Delta+1-Q^2)\Phi_{I}(x)=-2(x^{\bot}\cdot\Psi_{1})\ \ \hbox{in}\ \ \R^2,\ \int_{\R^2}Q\Phi_{I}dx=0,
\end{equation}
where $\Psi_{1}$ satisfies \eqref{0.1}.

Under the assumption  \eqref{0:1}, i.e., $\Omega_{a_{k}}=C_{0}(a^*-a_{k})^{-\beta}$ for some $\beta\in[0,\frac{1}{2})$ and $C_{0}>0$, we shall establish in Theorem \ref{th5}  the following refined expansion of the minimizer $u_{a_{k}}$ for $I(a_{k},\Omega_{a_{k}})$:
\begin{equation}\label{E1}
\begin{split}
\nu_{a_{k}}(x):&=\bar{\eps}_{a_{k}}\sqrt{a_{k}}u_{a_{k}}(\bar{\eps}_{a_{k}}x+x_{a_{k}})
e^{-i(\Omega_{a_{k}}\bar{\eps}_{a_{k}}x\cdot x_{a_{k}}^{\bot}-\rho_{a_{k}})}\\
&=Q(x)+\Omega_{a_{k}}^2\bar{\eps}_{a_{k}}^4\Psi_{1}(x)
+\Omega_{a_{k}}^2\bar{\eps}_{a_{k}}^5\Psi_{2}(x)\\
&\quad+i\Omega_{a_{k}}^3\bar{\eps}_{a_{k}}^6\Phi_{I}(x)+o(\Omega_{a_{k}}^3\bar{\eps}_{a_{k}}^6)
\ \ \hbox{as}\ \ k\to\infty,
\end{split}
\end{equation}
where  $\rho_{a_{k}}\in [0,2\pi)$ is a suitable constant phase, and
\begin{equation}\label{0:2}
    \bar{\eps}_{a_{k}}:=\sqrt{-\frac{1}{\mu_{a_{k}}}}=\frac{(a^*-a_{k})^{\frac{1+2\beta}{4}}}{\sqrt{C_{0}}\lambda}\big[1+o(1)\big]\ \ \hbox{as}\ \ k\to\infty.
\end{equation}
Moreover, the unique global maximum point $x_{a_{k}}$ of $|u_{a_{k}}|$ satisfies
 \begin{equation*}
    |x_{a_{k}}|^2-1=\tilde{C}\bar{\eps}_{a_{k}}^2+o(\Omega_{a_{k}}\bar{\eps}_{a_{k}}^4)\ \ \hbox{as}\ \ k\to\infty,
 \end{equation*}
where the constant $\tilde{C}=-\frac{8\lambda^4}{5a^*}<0$ is as before. We remark from  (\ref{E1}) that the rotational velocity $\Omega_{a_{k}}$ affects the minimizer $u_{a_{k}}$ starting from the second term of  the expansion (\ref{E1}). This is different from the harmonic cases discussed in \cite[Theorem 1.3]{GLY} and \cite[Theorem 1.1]{G}. Since the refined expansion (\ref{E1}) is established in terms of $\bar{\eps}_{a_{k}}$ defined by (\ref{0:2}), instead of $\eps_{a_{k}}:=\frac{1}{\sqrt{C_{0}}\lambda}(a^*-a_{k})^{\frac{1+2\beta}{4}}$, one does not need as in \cite{G} to further estimate the Lagrange multiplier $\mu_{a}$ for \eqref{ELE1}.
Comparing with \cite[Theorem 1.2]{G}, this simplifies greatly the whole proof of  Theorem \ref{th4}.

As the last step, we shall complete the proof of  Theorem \ref{th4} by applying \eqref{E1}. In fact, under the assumptions of Theorem \ref{th4}, the proof of \eqref{E1} implies that for sufficiently large constant $R>1$, \begin{equation}\label{AA:6.38}
|\nu_{a_{k}}(x)|\geq Q(x)-C_{R}\Omega_{a_{k}}^2\bar{\eps}_{a_{k}}^4>0\ \ \hbox{in}\ \ \big\{x\in\R^2:|x|\leq R\big\}\ \ \hbox{as}\ \ k\to\infty,
\end{equation}
where $\bar{\eps}_{a_{k}}$ is as in \eqref{E1} and  (\ref{0:2}).
On the other hand, the difference function
$$\tilde{w}_{a_{k}}(x):=\nu_{a_{k}}(x)-Q(x):
=\bar{\eps}_{a_{k}}\sqrt{a_{k}}u_{a_{k}}(\bar{\eps}_{a_{k}}x+x_{a_{k}})e^{-i(\bar{\eps}_{a_{k}}\Omega_{a_{k}} x\cdot x_{a_{k}}^{\bot}-\rho_{a_{k}})}-Q(x)$$
satisfies
\begin{equation}\label{1.06}
\begin{split}
-\Delta|\tilde{w}_{a_{k}}|+|\tilde{w}_{a_{k}}|
 \leq & \bar{C}_{0} \big(\Omega_{a_{k}}^2\bar{\eps}_{a_{k}}^4|x|^\frac{3}{2}
+\Omega_{a_{k}}^2\bar{\eps}_{a_{k}}^5|x|^{\frac{5}{2}}
+\Omega_{a_{k}}^2\bar{\eps}_{a_{k}}^6|x|^{\frac{7}{2}}\big)e^{-|x|}\\
    \ \ &\hbox{in}\ \ \R^2\backslash B_{R}(0)\ \ \hbox{as}\ \ k\to\infty,
\end{split}
\end{equation}
where $\bar{\eps}_{a_{k}}>0$ and $\rho_{a_{k}}\in [0,2\pi)$ are as above.

The proof of \eqref{E1} shows that for above sufficiently large constant $R>1$,
\begin{equation*}
  |\tilde{w}_{a_{k}}|\leq\bar{C} \big(\Omega_{a_{k}}^2\bar{\eps}_{a_{k}}^4|x|^\frac{5}{2}
+\Omega_{a_{k}}^2\bar{\eps}_{a_{k}}^5|x|^{\frac{7}{2}}
+\Omega_{a_{k}}^2\bar{\eps}_{a_{k}}^6|x|^{\frac{9}{2}}\big)e^{-|x|}\ \ \hbox{at}\ \ |x|=R>1\ \ \hbox{as}\ \ k\to\infty.
\end{equation*}
Different from \cite[Theorem 1.1]{G}, one however cannot estimate directly the right hand side of \eqref{1.06}. To overcome this difficulty, defining
\begin{equation*}
    \tilde{w}_{1a_{k}}=C_{1}\Omega_{a_{k}}^2\bar{\eps}_{a_{k}}^4|x|^\frac{5}{2}e^{-|x|},\ \, \tilde{w}_{2a_{k}}=C_{2}\Omega_{a_{k}}^2\bar{\eps}_{a_{k}}^5|x|^\frac{7}{2}e^{-|x|},\ \, \hbox{and}\, \ \tilde{w}_{3a_{k}}=C_{3}\Omega_{a_{k}}^2\bar{\eps}_{a_{k}}^6|x|^\frac{9}{2}e^{-|x|}\, \ \hbox{in}\, \ \R^2,
\end{equation*}
where $C_{i}>0$ are properly chosen constants, we shall prove in Section 5 that
\begin{equation*}
    |\tilde{w}_{a_{k}}|\leq \tilde{w}_{1a_{k}}+\tilde{w}_{2a_{k}}+\tilde{w}_{3a_{k}}\ \ \hbox{in}\ \ \R^2\backslash B_{R}(0)\ \ \hbox{as}\ \ k\to\infty.
\end{equation*}
By the definition of $\tilde{w}_{ia_{k}}$ and the exponential decay of $Q(x)$, we shall further derive that
\begin{equation*}
  \frac{1}{3}Q-\tilde{w}_{ia_{k}}>0,\ \ \hbox{if} \ \ R\leq|x|\leq C_{4}(\Omega_{a_{k}}^2\bar{\eps}_{a_{k}}^6)^{-\frac{1}{5}}\ \ \hbox{as}\ \  k\to\infty,\,\ i=1, 2, 3.
\end{equation*}
This hence implies that
\begin{equation*}
\begin{aligned}
  |\nu_{a_{k}}|&\geq Q-|\nu_{a_{k}}-Q|=Q-|\tilde{w}_{a_{k}}|
  \geq \sum_{i=1}^3\Big(\frac{1}{3}Q-\tilde{w}_{ia_{k}}\Big)>0,\\
  &\ \ \hbox{if}\ \ R\leq|x|\leq C_{4}(\Omega_{a_{k}}^2\bar{\eps}_{a_{k}}^6)^{-\frac{1}{5}}\ \ \hbox{as}\ \  k\to\infty.
\end{aligned}
\end{equation*}
Together with \eqref{AA:6.38}, this thus yields that $|\nu_{a_{k}}|>0$ holds for $|x|\leq C_{4}(\Omega_{a_{k}}^2\bar{\eps}_{a_{k}}^6)^{-\frac{1}{5}}$ as $k\to\infty$. We therefore complete the proof of Theorem \ref{th4} in view of \eqref{E1} and \eqref{0:2}.

This paper is organized as follows: In Section 2, we mainly establish Theorem \ref{th2.1} on the refined energy estimates of $I(a,\Omega_{a})$ as $a\nearrow a^*$. The proof of Theorem \ref{th1} is then given in Section 3 by employing energy methods and blow-up analysis. In Section 4, we shall derive the refined estimates of minimizers $u_{a}$ as $a\nearrow a^*$. Based on the estimates established in previous sections, in Section 5 we first prove Theorem \ref{th5} on the refined expansion of $u_{a}$ as $a\nearrow a^*$, by which we finally complete the proof of Theorem \ref{th4} on the non-existence of vortices in the larger region $R(a)$.

\section{Refined Energy Estimates of $I(a,\Omega_{a})$ as $a\nearrow a^*$}

In this section, we first establish the following Theorem \ref{th2.1} on the energy estimates of  $I(a,\Omega_{a})$ as $a\nearrow a^*$,  where $\Omega_{a}>0$ satisfies \eqref{0:1}. We then employ it to analyze the blow up rates and the maximum point of the minimizer $|u_{a}|$ for $I(a,\Omega_{a})$ as $a\nearrow a^*$.

\begin{thm}\label{th2.1}
There exist two constants $C_{1}>0$ and $C_{2}>0$, independent of $0< a<a^*$, such that
\begin{equation*}
    C_{1}(a^*-a)^{\frac{1}{2}-\beta}\leq I(a,\Omega_{a})\leq C_{2}(a^*-a)^{\frac{1}{2}-\beta}\ \ \hbox{as}\ \ a\nearrow a^*,
\end{equation*}
where $\Omega_{a}>0$ satisfies \eqref{0:1}, $i.e.,$ $\Omega_{a}:=C_{0}(a^*-a)^{-\beta}>0$ holds for some $\beta\in [0,\frac{1}{2})$ and $C_{0}>0$.
\end{thm}

To prove Theorem \ref{th2.1}, we recall from \cite[Lemma 8.1.2]{CAZ} that the unique positive solution $Q(x)$ of (\ref{Q}) satisfies
\begin{equation}\label{Q1}
    \int_{\R^2}|\nabla Q|^2dx=\int_{\R^2}Q^2dx=\frac{1}{2}\int_{\R^2}Q^4dx.
\end{equation}
Furthermore, we have the following exponential decay
\begin{equation}\label{expo}
    Q(x),\ \  |\nabla Q(x)|=O\big(|x|^{-\frac{1}{2}}e^{-|x|}\big)\ \ \hbox{as}\ \ |x|\rightarrow\infty.
\end{equation}
It was proved in \cite{W} that the identity of the following Gagliardo-Nirenberg inequality
\begin{equation}\label{GN}
    \int_{\R^2}|u|^{4}dx\leq\frac{2}{a^*}\int_{\R^2}|\nabla u|^2dx\int_{\R^2}|u|^{2}dx,\ \ u\in H^{1}(\R^2,\R)
\end{equation}
is achieved at $u(x)=Q(x)$. We now start with the following imprecise bounds of the energy $I(a,\Omega_{a})$ as $a\nearrow a^*$.

\begin{lem}\label{lem2.1}
There exist two constants $C_{1}'>0$ and $C_{2}'>0$, independent of $0< a<a^*$, such that
\begin{equation}\label{est}
   0< C_{1}'(a^*-a)^{\frac{2}{3}(1-\beta)}\leq I(a,\Omega_{a})\leq C_{2}'(a^*-a)^{\frac{1}{2}-\beta}\ \ \
    \hbox{as}\ \ a\nearrow a^*,
\end{equation}
where $\Omega_{a}>0$ satisfies \eqref{0:1}, $i.e.,$ $\Omega_{a}:=C_{0}(a^*-a)^{-\beta}>0$ holds for some $\beta\in [0,\frac{1}{2})$ and $C_{0}>0$.
\end{lem}

\noindent \textbf{Proof.}
For any $\sigma>0$, consider $u\in\mathcal{H}$ satisfying $\|u\|_{2}^2=1$. By the Gagliardo-Nirenberg inequality (\ref{GN}) and the diamagnetic inequality (\ref{dia}), we have
\begin{equation}\label{lest}
\begin{split}
     F_{a,\Omega_{a}}(u)
     &=\int_{\R^{2}}\big|(\nabla -i\Omega_{a} x^{\bot})u\big|^{2}dx+\frac{\Omega_{a}^{2}}{8}\int_{\R^{2}}(|x|^{2}-1)^{2}|u|^{2}dx-\frac{a}{2}\int_{\R^{2}}|u|^{4}dx\\
     &\geq\int_{\R^{2}}\big|\nabla|u|\big|^{2}dx+\frac{\Omega_{a}^{2}}{8}\int_{\R^{2}}(|x|^{2}-1)^{2}|u|^{2}dx-\frac{a}{2}\int_{\R^{2}}|u|^{4}dx\\
     &\geq \frac{\Omega_{a}^{2}}{8}\sigma-\frac{\Omega_{a}^{2}}{8}\int_{\R^{2}}\big[\sigma-(|x|^{2}-1)^{2}
     \big]_{+}|u|^{2}dx+\frac{a^*-a}{2}\int_{\R^{2}}|u|^{4}dx\\
     &\geq \frac{\Omega_{a}^{2}}{8}\sigma-\frac{\Omega_{a}^{4}}{128(a^*-a)}\int_{\R^{2}}\big[\sigma-(|x|^{2}-1)^{2}\big]^{2}_{+}dx,
\end{split}
\end{equation}
where $[\cdot]=max\{0,\cdot\}$ denotes the positive part. For $\sigma>0$ small enough, we have
\begin{equation}\label{est2}
\begin{split}
    \int_{\R^{2}}\big[\sigma-(|x|^{2}-1)^{2}\big]^{2}_{+}dx
    &=2\pi\int^{(1+\sigma^{\frac{1}{2}})^{\frac{1}{2}}}_{(1-\sigma^{\frac{1}{2}})^{\frac{1}{2}}}\big[\sigma-(r^{2}-1)^{2}\big]^{2}rdr\\
    &=\pi\int^{\frac{\pi}{2}}_{-\frac{\pi}{2}}\sigma^{2}\cos^{4}\theta
    \big(1+\sqrt{\sigma}
    \sin\theta\big)^{\frac{1}{2}}\big(1+\sqrt{\sigma}\sin\theta\big)^{-\frac{1}{2}}\sqrt{\sigma}\cos\theta d\theta\\
    &=\pi\sigma^{\frac{5}{2}}\int^{\frac{\pi}{2}}_{-\frac{\pi}{2}}\cos^{5}\theta d\theta
    = \frac{16\pi}{15}\sigma^{\frac{5}{2}},
    \end{split}
\end{equation}
where we set $r=\big(1+\sigma^{\frac{1}{2}}\sin\theta\big)^{\frac{1}{2}}$ and $-\frac{\pi}{2}\leq\theta\leq\frac{\pi}{2}$ in the second identity.
Combining (\ref{lest}) and (\ref{est2}), it then yields that
\begin{equation}\label{lest1}
\begin{split}
    F_{a,\Omega_{a}}(u)
    &\geq \frac{\Omega_{a}^{2}}{8}\sigma-\frac{\pi\Omega_{a}^{4}}{120(a^*-a)}\sigma^{\frac{5}{2}}
    \geq C_{1}'(a^*-a)^{\frac{2}{3}(1-\beta)},
\end{split}
\end{equation}
by taking $\sigma=\big[\frac{6(a^*-a)}{\pi\Omega_{a}^{2}}\big]^{\frac{2}{3}}$, and the lower bound of (\ref{est}) is therefore proved.

We next establish the upper bound of (\ref{est}).
Choose a non-negative function $\varphi\in C^{\infty}_{0}(\R^2)$ satisfying $\varphi=1$ for $|x|\leq1$ and $\varphi=0$ for $|x|\geq2$. For any $\tau>0$, we define
\begin{equation*}
   u_{\tau}(x)=A_{\tau}\frac{\tau}{\|Q\|_{2}}\varphi(x-y_{0})Q\big(\tau (x-y_{0})\big)e^{i\Omega_{a} S(x)},
\end{equation*}
where $S(x)=x\cdot y^{\bot}_{0}$, $|y_{0}|=1$, and $A_{\tau}>0$ is chosen such that $\int_{\R^2}|u_{\tau}|^{2}dx=1$. Then $A_{\tau}$ satisfies
\begin{equation*}
    \frac{1}{A_{\tau}^{2}}=\frac{1}{\|Q\|_{2}^{2}}\int_{\R^2}Q^2(x)\varphi^{2}\big(\frac{x}{\tau}\big)dx=1+O(\tau^{-\infty})\ \ \hbox{as}\ \ \tau\rightarrow\infty,
\end{equation*}
due to the exponential decay (\ref{expo}) of $Q$. Here and below we use $f(t)=O(t^{-\infty})$ to denote a function $f$ satisfying $\lim\limits_{t\rightarrow\infty}|f(t)|t^s=0$ for all $s>0$. We then get that
\begin{equation*}
\begin{aligned}
    F_{a,\Omega_{a}}(u_{\tau})&=\int_{\R^2}\big|(\nabla-i\Omega_{a} x^{\bot})u_{\tau}\big|^2dx+\frac{\Omega_{a}^2}{8}\int_{\R^{2}}\big(|x|^2-1\big)^2|u_{\tau}|^2dx-\frac{a}{2}\int_{\R^2}|u_{\tau}|^{4}dx\\
    &=\frac{\tau^{2}}{a^*}\int_{\R^2}\Big|\nabla\Big(Q\big(\tau(x-y_{0})\big)e^{i\Omega_{a} S(x)}\Big)-i\Omega_{a} x^{\bot}Q\big(\tau(x-y_{0})\big)e^{i\Omega_{a} S(x)}\Big|^{2}dx\\
   &\quad+\frac{\Omega_{a}^{2}\tau^2}{8a^*}\int_{\R^2}\big(|x|^2-1\big)^{2}Q^{2}\big(\tau (x-y_{0})\big)dx-\frac{\tau^{4}a}{2(a^{*})^2}\int_{\R^2}Q^{4}\big(\tau(x-y_{0})\big)dx\\
   &\quad+O(\tau^{-\infty})\\
   &=\frac{\tau^{2}}{a^*}\int_{\R^2}\Big|\nabla Q
  -\frac{i\Omega_{a}}{\tau^{2}}x^{\bot}Q(x)\Big|^{2}dx
  +\frac{\Omega_{a}^{2}}{8a^*}\int_{\R^2}\Big(\big|\frac{x}{\tau}+y_{0}\big|^2-1\Big)^{2}Q^2(x)dx
  -\frac{\tau^{2}a}{a^*}\\
  &\quad+O(\tau^{-\infty})\\
  &=\frac{\tau^{2}(a^*-a)}{a^*}+\frac{\Omega_{a}^{2}}{\tau^{2}a^*}\int_{\R^2}|x|^{2}Q^{2}dx
  +\frac{\Omega_{a}^{2}}{8a^*}\int_{\R^2}\Big(\big|\frac{x}{\tau}+y_{0}\big|^2-1\Big)^{2}Q^2(x)dx\\
  &\quad+O(\tau^{-\infty})\ \ \mbox{as}\ \ \tau\rightarrow\infty.
\end{aligned}
\end{equation*}
This gives that
\begin{equation}\label{higest1}
\begin{split}
    I(a,\Omega_{a})&\leq F_{a,\Omega_{a}}(u_{\tau})\leq C(a^*-a)\tau^2+C^{'}\Omega_{a}^{2}\tau^{-2}
    \leq C_{2}'(a^*-a)^{\frac{1}{2}-\beta},\\
\end{split}
\end{equation}
by taking $\tau=(a^*-a)^{-\frac{2\beta+1}{4}}$.
Applying (\ref{lest1}) and (\ref{higest1}), the proof of Lemma \ref{lem2.1} is therefore complete.
\qed
\vskip 0.05truein

To complete the proof of Theorem \ref{th2.1}, one can note from Lemma \ref{lem2.1} that it is still necessary to derive the precise lower bound of $I(a,\Omega_{a})$ as $ a \nearrow a^*$. For this reason, we note that if $u_{a}$ is a minimizer of $I(a,\Omega_{a})$ for $0< a<a^*$, then $u_{a}$ satisfies the elliptic equation (\ref{ELE1}), where $\mu_{a}\in\R$ is the associated Lagrange multiplier. We further deduce from (\ref{ELE1}) that $\mu_{a}$ satisfies
\begin{equation}\label{mu}
    \mu_{a}=I(a,\Omega_{a})-\frac{a}{2}\int_{\R^2}|u_{a}|^{4}dx,\ \ 0< a<a^*.
\end{equation}
We next analyze some properties of the minimizer $u_{a}$ for $I(a,\Omega_{a})$ with $0<a<a^*$.

\begin{lem}\label{lem2.3}
Let $u_a$ be a complex-valued minimizer of $I(a,\Omega_{a})$ for $0<a<a^*$, where $\Omega_{a}>0$ satisfies (\ref{0:1}) for some $\beta\in[0,\frac{1}{2})$ and $C_{0}>0$. Set
\begin{equation}\label{ypx7}
\epsilon_a:=\Big(\int_{\R^2}\big|\nabla |u_a| \big|^2 dx\Big)^{-\frac{1}{2}}>0.
\end{equation}
Then we have
\begin{enumerate}
\item The constant $\epsilon_a>0$ satisfies
  \begin{equation}\label{2:6}
    \epsilon_{a}\rightarrow 0\ \ \mbox{and}\ \ \mu_{a}\epsilon_{a}^{2}\rightarrow-1\ \ \mbox{as}\ \ a\nearrow a^*.
  \end{equation}
\item Define
  \begin{equation}\label{wa1}
    w_{a}(x):=\epsilon_{a}u_{a}(\epsilon_{a}x+x_{a})e^{-i(\epsilon_{a}\Omega_{a} x\cdot x_{a}^{\bot}-\theta_{a})},
  \end{equation}
where $x_{a}$ is a global maximum point of $|u_{a}|$ as $ a\nearrow a^*$, and $\theta_{a}\in [0,2\pi)$ is a proper constant. Then there exists a constant $\eta>0$, independent of $0<a<a^*$, such that
\begin{equation}\label{non01}
    \int_{B_{2}(0)}|w_{a}(x)|^2dx\geq\eta>0\ \ \hbox{as}\ \ a\nearrow a^*.
  \end{equation}

\item Any global maximum point $x_{a}$ of $|u_{a}|$ satisfies $\lim\limits_{a\nearrow a^*}|x_{a}|=1$. Moreover, for any sequence $\{a_{k}\}$ with $a_{k}\nearrow a^*$ as $k\to\infty$, there exists a subsequence, still denoted by $\{a_{k}\}$, of $\{a_{k}\}$ such that
      \begin{equation}\label{3:1}
        \lim\limits_{k\to\infty}x_{a_{k}}=x_{0}
      \end{equation}
      for some point $x_{0}\in\R^2$ satisfying $|x_{0}|=1$,
      and $w_{a_{k}}(x)$ satisfies
      \begin{equation}\label{conv7}
        \lim\limits_{k\to\infty} |w_{a_{k}}(x)|=w_{0}:=\frac{Q(x)}{\sqrt{a^*}}\ \ \hbox{strongly in}\ \ H^{1}(\R^2).
      \end{equation}
\end{enumerate}
\end{lem}

\noindent \textbf{Proof.}
1. By contradiction, assume that $\epsilon_{a}\nrightarrow 0$ as $a\nearrow a^*$. The Gagliardo-Nirenberg inequality (\ref{GN}) and the diamagnetic inequality (\ref{dia}) yield that
\begin{equation}\label{ene}
\begin{split}
    I(a,\Omega_{a})&=\int_{\R^2}\big|(\nabla-i\Omega_{a} x^{\bot})u_{a}\big|^{2}dx
    +\frac{\Omega_{a}^{2}}{8}\int_{\R^2}(|x|^{2}-1)^{2}|u_{a}|^{2}dx
    -\frac{a}{2}\int_{\R^2}|u_{a}|^{4}dx\\
    &\geq \frac{\Omega_{a}^{2}}{8}\int_{\R^2}(|x|^{2}-1)^{2}|u_{a}|^{2}dx.
\end{split}
\end{equation}
Applying Lemma \ref{lem2.1}, we then deduce from (\ref{ene}) that $\{|u_{a}|\}$ is bounded uniformly in $\mathcal{H}$ as $a\nearrow a^*$. Since the embedding $\mathcal{H}\hookrightarrow L^{q}(\R^2,\C)$ is compact for $2\leq q<+\infty$ (similar to \cite[Lemma 3.1]{ZJ}), there exist a subsequence, still denoted by $\{|u_{a}|\}$, of $\{|u_{a}|\}$ and a function $0\leq u_{0}\in \mathcal{H}$ such that
\begin{equation*}
    |u_{a}|\rightharpoonup u_{0}\ \ \hbox{weakly in}\ \ \mathcal{H}\ \ \hbox{and}\ \ |u_{a}|\rightarrow u_{0}\ \ \hbox{strongly in}\ \ L^{q}(\R^2)\ \ \hbox{as}\ \ a\nearrow a^*,
\end{equation*}
where $q\in[2,+\infty)$. Since $|u_{a}|\rightarrow u_{0}$ strongly in $L^{2}(\R^2)$ as $a\nearrow a^*$, we obtain that $\|u_{0}\|_{L^2(\R^2)}=1$. Following (\ref{est}) and (\ref{ene}), we have
\begin{equation*}
    C(a^*-a)^{\frac{1}{2}-\beta}\geq I(a,\Omega_{a})\geq \frac{\Omega_{a}^{2}}{8}\int_{\R^2}(|x|^{2}-1)^{2}|u_{a}|^{2}dx,
\end{equation*}
which gives that
\begin{equation*}
    0=\lim\limits_{a\nearrow a^*}\frac{8I(a,\Omega_{a})}{\Omega_{a}^{2}}\geq\liminf_{a\nearrow a^*}\int_{\R^2}(|x|^{2}-1)^{2}|u_{a}|^{2}dx
    \geq\int_{\R^2}(|x|^{2}-1)^{2}|u_{0}|^{2}dx\geq 0.
\end{equation*}
This implies that $u_{0}\equiv 0$ in $\R^2$, a contradiction. We therefore have $\epsilon_{a}\rightarrow 0$ as $a\nearrow a^*$.

In the following, we prove that $\mu_{a}\epsilon_{a}^{2}\rightarrow -1$ as $a\nearrow a^*$.
Since $\epsilon_{a}\rightarrow 0$ as $a\nearrow a^*$, we obtain from the  inequality (\ref{GN}) and Lemma \ref{lem2.1} that
\begin{equation}\label{-1}
    \limsup_{a\nearrow a^*}\frac{\frac{a}{2}\int_{\R^2}|u_{a}|^{4}dx-I(a,\Omega_{a})}{\int_{\R^2}\big|\nabla |u_{a}|\big|^{2}dx}\leq 1.
\end{equation}
Assume that $\mu_{a}\epsilon_{a}^{2}\nrightarrow -1$ as $a\nearrow a^*$. Up to a subsequence if necessary, we then obtain from (\ref{mu}) and (\ref{-1}) that there exists a constant $\gamma_{0}>0$ such that
\begin{equation*}
    -\liminf_{a\nearrow a^*}\mu_{a}\epsilon_{a}^{2}=\limsup_{a\nearrow a^*}\frac{\frac{a}{2}\int_{\R^2}|u_{a}|^{4}dx-I(a,\Omega_{a})}{\int_{\R^2}\big|\nabla |u_{a}|\big|^{2}dx}\leq 1-\gamma_{0},
\end{equation*}
from which we thus deduce that
\begin{equation}\label{contra}
\begin{split}
   I(a,\Omega_{a})&=\int_{\R^2}\big|(\nabla-i\Omega_{a} x^{\bot})u_{a}\big|^{2}dx
    +\frac{\Omega_{a}^{2}}{8}\int_{\R^2}(|x|^{2}-1)^{2}|u_{a}|^{2}dx
    -\frac{a}{2}\int_{\R^2}|u_{a}|^{4}dx\\
    &\geq\int_{\R^2}\big|\nabla |u_{a}|\big|^{2}dx
    -\frac{a}{2}\int_{\R^2}|u_{a}|^{4}dx
    \geq \gamma_{0}\int_{\R^2}\big|\nabla |u_{a}|\big|^{2}dx.
\end{split}
\end{equation}
Since $\int_{\R^2}\big|\nabla |u_{a}|\big|^{2}dx\rightarrow\infty$ as $a\nearrow a^*$, we deduce from (\ref{contra}) that
$I(a,\Omega_{a})\rightarrow\infty$ as $a\nearrow a^*$, which however contradicts with Lemma \ref{lem2.1}. Hence, $\mu_{a}\epsilon_{a}^2\rightarrow -1$ as $a\nearrow a^*$.

2. Denote $\bar{w}_{a}(x):=\epsilon_{a}u_{a}(\epsilon_{a}x+x_{a})e^{-i(\epsilon_{a}\Omega_{a} x\cdot x_{a}^{\bot})}$ and
$w_{a}(x):=\bar{w}_{a}(x)e^{i\theta_{a}}$, where the parameter $\theta_{a}\in [0,2\pi)$ is chosen such that
\begin{equation}\label{orth3}
    \Big\|w_{a}-\frac{Q}{\sqrt{a^*}}\Big\|_{L^{2}(\R^2)}=\min\limits_{\theta\in [0,2\pi)}\Big\|e^{i\theta}\bar{w}_{a}-\frac{Q}{\sqrt{a^*}}\Big\|_{L^{2}(\R^2)}.
\end{equation}
Rewrite $w_{a}(x):=R_{a}(x)+iI_{a}(x)$, where $R_{a}(x)$ denotes the real part of $w_{a}(x)$, and $I_{a}(x)$ denotes the imaginary part of $w_{a}(x)$.
From \eqref{orth3}, we then obtain the following orthogonality condition:
\begin{equation}\label{orth4}
    \int_{\R^2}Q(x)I_{a}(x)dx=0.
\end{equation}
It follows from (\ref{ELE1}) and \eqref{wa1} that $w_{a}(x)$ satisfies the following equation
\begin{equation}\label{waele}
    -\Delta w_{a}+2i\Omega_{a}\epsilon_{a}^{2}(x^{\bot}\cdot\nabla w_{a})+\Big[\Omega_{a}^{2}\epsilon_{a}^{4}|x|^{2}+\frac{\Omega_{a}^{2}\epsilon_{a}^{2}}{8}\big(|\epsilon_{a}x+x_{a}|^{2}-1\big)^{2}
    -\mu_{a}\epsilon_{a}^{2}-a|w_{a}|^{2}\Big]w_{a}=0\ \ \hbox{in}\ \ \R^2.
\end{equation}

Defining $W_{a}(x):=|w_{a}(x)|^{2}\geq 0$ in $\R^2$, we then derive from (\ref{waele}) that
\begin{equation}\label{Wa}
    \begin{split}
    &-\frac{1}{2}\Delta W_{a}(x)+|\nabla w_{a}|^{2}-2\Omega_{a}\epsilon_{a}^{2}x^{\bot}\cdot(iw_{a},\nabla w_{a})\\
    &+\Big[\Omega_{a}^{2}\epsilon_{a}^{4}|x|^{2}
    +\frac{\Omega_{a}^{2}\epsilon_{a}^{2}}{8}\big(|\epsilon_{a}x+x_{a}|^{2}-1\big)^{2}-\mu_{a}\epsilon_{a}^{2}-aW_{a}\Big]W_{a}=0
    \ \ \hbox{in}\ \ \R^2.
\end{split}
\end{equation}
By the diamagnetic inequality (\ref{dia}), we have
\begin{equation*}
    |\nabla w_{a}|^{2}-2\Omega_{a}\epsilon_{a}^{2}x^{\bot}\cdot(iw_{a},\nabla w_{a})+\Omega_{a}^{2}\epsilon_{a}^{4}|x|^{2}W_{a}\geq 0 \ \ \hbox{in}\ \ \R^2,
\end{equation*}
together with \eqref{Wa}, which yields that
\begin{equation}\label{Wa4}
    -\frac{1}{2}\Delta W_{a}(x)-\mu_{a}\epsilon_{a}^{2}W_{a}-aW_{a}^{2}\leq 0 \ \ \hbox{in}\ \ \R^2.
\end{equation}
Since the origin is a global maximum point of $W_{a}(x)$ for all $0< a<a^*$, we have $-\Delta W_{a}(0)\geq0$ for all $0< a<a^*$. By the fact that $\mu_{a}\epsilon_{a}^{2}\rightarrow -1$ as $a\nearrow a^*$, we then get from \eqref{Wa4} that $W_{a}(0)\geq \gamma>0$ holds uniformly as $a\nearrow a^*$. Following De Giorgi-Nash-Moser theory \cite[Theorem 4.1]{HAN}, we obtain from (\ref{Wa4}) that
\begin{equation}\label{Wa3}
    \int_{B_{2}(0)}W_{a}^{2}dx\geq C\max\limits_{x\in B_{1}(0)}W_{a}(x)\geq C_{1}(\gamma)\ \ \ \hbox{as}\ \ a\nearrow a^*.
\end{equation}
Moreover, since $\{|w_{a}|\}$ is bounded uniformly in $H^{1}(\R^2)$, we derive from (\ref{Wa3}) that
\begin{equation}\label{Wa0}
    W_{a}(0)=\max\limits_{x\in\R^2} W_{a}(x)\leq C\int_{B_{2}(0)}W_{a}^{2}dx\leq C\ \ \hbox{as}\ \ a\nearrow a^*.
\end{equation}
We conclude from (\ref{Wa3}) and $(\ref{Wa0})$ that
\begin{equation*}
    \int_{B_{2}(0)}|w_{a}|^{2}dx\geq\int_{B_{2}(0)}\frac{W_{a}^{2}}{\max\limits_{x\in\R^2}W_{a}}dx\geq C_{2}(\gamma):=\eta>0\ \ \hbox{as}\ \ a\nearrow a^*,
\end{equation*}
which completes the proof of (\ref{non01}).

3.
We first prove that any global maximum point $x_{a}$ of $|u_{a}|$ satisfies
\begin{equation}\label{2:4}
    \lim\limits_{a\nearrow a^*}|x_{a}|=1.
\end{equation}
Indeed, if \eqref{2:4} is false, then there exist a constant $\alpha>0$ and a subsequence $\{a_{k}\}$, where $a_{k}\nearrow a^*$ as $k\rightarrow\infty$, such that
\begin{equation*}
    \big||x_{a_{k}}|^2-1\big|\geq\alpha>0\ \ \hbox{as}\ \ k\to\infty.
\end{equation*}
By Fatou's lemma, it therefore follows from \eqref{2:6} and (\ref{non01}) that
\begin{equation}\label{2:5}
\begin{aligned}
    \liminf_{k\rightarrow\infty}\int_{\R^2}\big||\epsilon_{a_{k}}x+x_{a_{k}}|^2-1\big|^2|w_{a_{k}}|^2dx
    &\geq\int_{B_{2}(0)}\liminf_{k\rightarrow\infty}\big||\epsilon_{a_{k}}x+x_{a_{k}}|^2-1\big|^2|w_{a_{k}}|^2dx\\
    &\geq \frac{\alpha^2\eta}{2}>0.
\end{aligned}
\end{equation}
However, we deduce from (\ref{ene}) and Lemma \ref{lem2.1} that
\begin{equation*}
    \lim\limits_{k\rightarrow\infty}\int_{\R^2}(|x|^2-1)^2|u_{a_{k}}|^2dx
    =\lim\limits_{k\rightarrow\infty}\int_{\R^2}\big||\epsilon_{a_{k}}x+x_{a_{k}}|^2-1\big|^2|w_{a_{k}}|^2dx=0,
\end{equation*}
which contradicts with (\ref{2:5}). Hence, the proof of \eqref{2:4} is complete, which further implies that \eqref{3:1} holds true.

Let $\{a_{k}\}$ be the convergent subsequence obtained in \eqref{3:1}. Since $\{|w_{a_{k}}|\}$ is bounded uniformly in $H^{1}(\R^2)$, there exists a subsequence, still denoted by $\{|w_{a_{k}}|\}$, of $\{|w_{a_{k}}|\}$ such that $|w_{a_{k}}|\rightharpoonup w_{0}$ weakly in $H^{1}(\R^2)$ as $k\to\infty$ for some $0\leq w_{0}\in H^{1}(\R^2)$. By (\ref{non01}), we may assume that $|w_{a_{k}}|\rightarrow w_{0}\not\equiv 0$ almost everywhere in $\R^2$. By Br\'ezis-Lieb lemma \cite{Willem}, we have
\begin{equation}\label{BL}
    \|w_{a_{k}}\|^{q}_{q}=\|w_{0}\|^{q}_{q}+\big\||w_{a_{k}}|-w_{0}\big\|^{q}_{q}+o(1)\ \ \hbox{as}\ \ k\to\infty,
    \ \ \hbox{where}\ \ 2\leq q<\infty,
\end{equation}
and
\begin{equation}\label{BL1}
    \big\|\nabla|w_{a_{k}}|\big\|^{2}_{2}=\|\nabla w_{0}\|^{2}_{2}+\big\|\nabla(|w_{a_{k}}|-w_{0})\big\|^{2}_{2}+o(1)\ \ \hbox{as}\ \ k\to\infty,
\end{equation}
which imply that $\big\||w_{a_{k}}|-w_{0}\big\|^{2}_{2}\leq 1$ as $k\to\infty$.
We obtain from Lemma \ref{lem2.1} that $I(a_{k},\Omega_{a_{k}})\rightarrow 0$ as $k\to\infty$, which implies from \eqref{2:6} that
\begin{equation*}
    1=-\lim\limits_{k\to\infty}\mu_{a_{k}}\epsilon_{a_{k}}^{2}
    =\lim\limits_{k\to\infty}\frac{\frac{a_{k}}{2}\int_{\R^2}|u_{a_{k}}|^{4}dx-I(a_{k},\Omega_{a_{k}})}{\int_{\R^2}\big|\nabla |u_{a_{k}}|\big|^{2}dx}.
\end{equation*}
By the definition of $w_{a_{k}}(x)$, we then have
\begin{equation*}
    \lim\limits_{k\to\infty}\int_{\R^2}|w_{a_{k}}|^{4}dx=\frac{2}{a^*},
\end{equation*}
which implies that
\begin{equation}\label{lim}
\begin{split}
    &\lim\limits_{k\to\infty}\Big\{\int_{\R^2}\big|\nabla|w_{a_{k}}|\big|^{2}dx
    -\frac{a^*}{2}\int_{\R^2}|w_{a_{k}}|^{4}dx\Big\}=0.
\end{split}
\end{equation}
Applying the Gagliardo-Nirenberg inequality (\ref{GN}), we deduce from (\ref{BL})--(\ref{lim}) that
\begin{equation}\label{lim1}
\begin{split}
    0&=\lim\limits_{k\to\infty}\Big\{\int_{\R^2}\big|\nabla|w_{a_{k}}|\big|^{2}dx
    -\frac{a^*}{2}\int_{\R^2}|w_{a_{k}}|^{4}dx\Big\}\\
    &=\int_{\R^2}|\nabla w_{0}|^2dx-\frac{a^*}{2}\int_{\R^2}|w_{0}|^4dx\\
    &\quad+\lim\limits_{k\to\infty}\Big\{\int_{\R^2}\big|\nabla(|w_{a_{k}}|-w_{0})\big|^{2}dx
    -\frac{a^*}{2}\int_{\R^2}\big||w_{a_{k}}|-w_{0}\big|^{4}dx\Big\}\\
    &\geq\int_{\R^2}|\nabla w_{0}|^2dx-\frac{a^*}{2}\int_{\R^2}|w_{0}|^4dx\\
    &\quad+\frac{a^*}{2}\lim\limits_{k\to\infty}\Big(\big\||w_{a_{k}}|-w_{0}\big\|^{-2}_{2}-1\Big)
    \int_{\R^2}\big||w_{a_{k}}|-w_{0}\big|^{4}dx\geq0.
\end{split}
\end{equation}
Using (\ref{GN}) again, it implies from (\ref{lim1}) that $\|w_{0}\|^{2}_{2}=1$, and hence
\begin{equation}\label{conv}
    |w_{a_{k}}|\rightarrow w_{0}\ \ \hbox{strongly in}\ \ L^2(\R^2)\ \ \hbox{as}\ \ k\to\infty.
\end{equation}
Since $\{|w_{a_{k}}|\}$ is bounded uniformly in $H^{1}(\R^2)$, we obtain from \eqref{conv} that $|w_{a_{k}}|\rightarrow w_{0}$ strongly in $L^4(\R^2)$
as $k\to\infty$.  It then follows from  (\ref{GN}) and (\ref{lim1}) that
 \begin{equation*}
    \lim\limits_{k\to\infty}\int_{\R^2}\big|\nabla |w_{a_{k}}|\big|^2= \frac{a^*}{2}\lim\limits_{k\to\infty}\int_{\R^2}|w_{a_{k}}|^4=\frac{a^*}{2}\int_{\R^2}|w_{0}|^4
    =\int_{\R^2}|\nabla w_{0}|^2,
 \end{equation*}
which implies that $\nabla |w_{a_{k}}|\rightarrow\nabla w_{0}$ strongly in $L^2(\R^2)$ as $k\to\infty$, and
$w_{0}$ must be an optimizer of the Gagliardo-Nirenberg inequality (\ref{GN}). Thus, there exists a point $z_{0}\in\R^2$ such that up to a subsequence if necessary,
\begin{equation*}
    |w_{a_{k}}(x)|\rightarrow w_{0}(x)=\frac{Q\big(|x+z_{0}|\big)}{\sqrt{a^*}}\ \ \hbox{strongly in}\ \ H^{1}(\R^2)\ \ \hbox{as}\ \ k\to\infty.
\end{equation*}

Since the origin is a global maximum point of $|w_{a_{k}}|$ for all $k>0$, the origin must be a global maximum point of $Q(x+z_{0})$, which implies that $z_{0}=0$, and up to a subsequence if necessary,
\begin{equation*}
    |w_{a_{k}}(x)|\rightarrow \frac{Q(x)}{\sqrt{a^*}}\ \ \hbox{strongly in}\ \ H^{1}(\R^2)\ \ \hbox{as}\ \ k\to\infty.
\end{equation*}
This completes the proof of Lemma \ref{lem2.3}.
\qed

\begin{lem}\label{lem2.6}
Under the assumptions of Lemma \ref{lem2.3}, let $\{a_{k}\}$ be the convergent subsequence obtained in Lemma \ref{lem2.3} (3).
Then there exists a constant $C>0$, independent of $a_{k}$, such that
\begin{equation*}
   \lim\limits_{k\to\infty}\frac{1}{\epsilon_{a_{k}}^{2}}\int_{\R^2}\big(|\epsilon_{a_{k}}x
    +x_{a_{k}}|^{2}-1\big)^{2}|w_{a_{k}}|^{2}dx\geq C,
\end{equation*}
where $\epsilon_{a_{k}}>0$ and $w_{a_{k}}$ are defined by \eqref{ypx7} and \eqref{wa1}, respectively.
\end{lem}

\noindent\textbf{Proof.}
For any $x\in\R^2$, let arg $x$ be the angle between $x$ and the positive $x$-axis, and denote $\langle x,y\rangle$ the angle between the vectors $x$ and $y$. We deduce from Lemma \ref{lem2.3} (3) that$\lim\limits_{k\to\infty}x_{a_{k}}=x_{0}$ and $|x_{0}|=1$. Without loss of generality, we may assume $x_{0}=(1,0)$, and hence arg $x_{a_{k}}\rightarrow0$ as $k\to\infty$. Choose $0<\delta<\frac{\pi}{16}$ small enough such that
\begin{equation}\label{arg}
    -\delta<\hbox{arg}\ x_{a_{k}}<\delta\ \ \hbox{as}\ \ k\to\infty.
\end{equation}
For any constant $R>0$, we denote
\begin{equation*}
\begin{split}
    \Omega_{\epsilon_{a_{k}}}^{1}&=\Big\{x\in B_{R}(0):\ |x|^2+\big|\frac{x_{a_{k}}}{\epsilon_{a_{k}}}\big|^2\leq\frac{1}{\epsilon_{a_{k}}^{2}}\Big\}\\
    &=\Big\{x\in B_{R}(0):\ |x|^2\leq\frac{1}{\epsilon_{a_{k}}^{2}}-\big|\frac{x_{a_{k}}}{\epsilon_{a_{k}}}\big|^2\Big\},\\
\end{split}
\end{equation*}
and
\begin{equation*}
\begin{split}
    \Omega_{\epsilon_{a_{k}}}^{2}&=\Big\{x\in B_{R}(0):\ |x|^2+\big|\frac{x_{a_{k}}}{\epsilon_{a_{k}}}\big|^2>\frac{1}{\epsilon_{a_{k}}^{2}}\Big\}\\
    &=\Big\{x\in B_{R}(0):\ \frac{1}{\epsilon_{a_{k}}^{2}}-\big|\frac{x_{a_{k}}}{\epsilon_{a_{k}}}\big|^2<|x|^2<R^2\Big\},\\
\end{split}
\end{equation*}
such that $B_{R}(0)=\Omega_{\epsilon_{a_{k}}}^{1}\cup\Omega_{\epsilon_{a_{k}}}^{2}$ and $\Omega_{\epsilon_{a_{k}}}^{1}\cap\Omega_{\epsilon_{a_{k}}}^{2}=\emptyset$. Since
\begin{equation*}
    |\Omega_{\epsilon_{a_{k}}}^{1}|+|\Omega_{\epsilon_{a_{k}}}^{2}|=|B_{R}(0)|=\pi R^2,
\end{equation*}
there exists a subsequence, still denoted by $\{\epsilon_{a_{k}}\}$, of $\{\epsilon_{a_{k}}\}$ such that
\begin{equation*}
    \hbox{either}\ \ |\Omega_{\epsilon_{a_{k}}}^{1}|\geq\frac{\pi R^2}{2}\ \
    \hbox{or}\ \ |\Omega_{\epsilon_{a_{k}}}^{2}|\geq\frac{\pi R^2}{2}.
\end{equation*}
We next consider the following two cases:

\emph{Case 1:} $|\Omega_{\epsilon_{a_{k}}}^{1}|\geq\frac{\pi R^2}{2}$. In this case, we have $B_{\frac{R}{\sqrt{2}}}(0)\subset\Omega_{\epsilon_{a_{k}}}^{1}$, and set
\begin{equation*}
    \Omega_{1}:=\Big(B_{\frac{R}{\sqrt{2}}}(0)\backslash B_{\frac{R}{2}}(0)\Big)
    \cap\Big\{x:\ \frac{\pi}{2}+2\delta<\hbox{arg}\ x<\frac{3\pi}{2}-2\delta\Big\}
    \subset\Omega_{\epsilon_{a_{k}}}^{1},
\end{equation*}
so that
\begin{equation*}
    |\Omega_{1}|=\frac{\pi-4\delta}{8}R^2.
\end{equation*}
From (\ref{arg}), we obtain that for any $x\in\Omega_{1}$,
\begin{equation*}
    x\cdot x_{a_{k}}=|x||x_{a_{k}}|\cos\langle x,x_{a_{k}}\rangle<0\ \
    \hbox{and}\ \ |\cos\langle x,x_{a_{k}}\rangle|>-\cos\Big(\frac{\pi}{2}+\delta\Big)>0.
\end{equation*}
Then for any $x\in\Omega_{1}$, there exists a constant $C>0$, independent of $a_{k}$, such that
\begin{equation*}
\begin{split}
    \frac{1}{\epsilon_{a_{k}}^{2}}\Big(|\epsilon_{a_{k}}x+x_{a_{k}}|^{2}-1\Big)^{2}
    &=\epsilon_{a_{k}}^{2}\Big[|x|^{2}+\big|\frac{x_{a_{k}}}{\epsilon_{a_{k}}}\big|^{2}+2|x| \big|\frac{x_{a_{k}}}{\epsilon_{a_{k}}}\big|\cos\langle x,x_{a_{k}}\rangle
    -\frac{1}{\epsilon_{a_{k}}^{2}}\Big]^{2}\\
    &\geq\epsilon_{a_{k}}^{2}\Big(2|x|\big|\frac{x_{a_{k}}}{\epsilon_{a_{k}}}\big|\cos\langle x,x_{a_{k}}\rangle\Big)^{2}\\
    &\geq4|x_{a_{k}}|^{2}|x|^{2}\cos^{2}\big(\frac{\pi}{2}+\delta\big)\\
    &\geq C|x|^2\cos^{2}\big(\frac{\pi}{2}+\delta\big)\ \ \hbox{as}\ \ k\rightarrow\infty,
\end{split}
\end{equation*}
due to the fact that $\lim\limits_{k\rightarrow\infty}|x_{a_{k}}|= 1$. Taking $\delta=\frac{\pi}{20}$, the above estimate thus implies that
\begin{equation}\label{2.01}
\begin{split}
    &\quad\lim\limits_{k\rightarrow\infty}\frac{1}{\epsilon_{a_{k}}^{2}}\int_{B_{R}(0)}\Big(|\epsilon_{a_{k}}x+x_{a_{k}}|^{2}-1\Big)^{2}|w_{a_{k}}|^{2}dx\\
    &\geq C \cos^{2}\frac{11\pi}{20}\int_{\Omega_{1}}|x|^{2}|w_{0}|^{2}dx
    \geq C(R),
\end{split}
\end{equation}
where $w_{0}>0$ is as in \eqref{conv7}.

\emph{Case 2:} $|\Omega_{\epsilon_{a_{k}}}^{2}|\geq\frac{\pi R^2}{2}$. In this case, we have the annular region $D_{R}:=B_{R}(0)\backslash B\frac{R}{\sqrt{2}}(0)\subset\Omega_{\epsilon_{a_{k}}}^{2}$. Set
\begin{equation*}
    \Omega_{2}:=D_{R}\cap \Big\{x:-\frac{\pi}{2}+2\delta<\hbox{arg}\ x<\frac{\pi}{2}-2\delta\Big\}\subset\Omega_{\epsilon_{a_{k}}}^{2},
\end{equation*}
so that
\begin{equation*}
    |\Omega_{2}|=\frac{\pi-4\delta}{4}R^2.
\end{equation*}
One can check that for any $x\in\Omega_{2}$,
\begin{equation*}
    x\cdot x_{a_{k}}=|x||x_{a_{k}}|\cos\langle x,x_{a_{k}}\rangle>0\ \
    \hbox{and}\ \ \cos\langle x,x_{a_{k}}\rangle>\cos\big(\frac{\pi}{2}-\delta\big)>0.
\end{equation*}
Then there exists a constant $C>0$, independent of $a_{k}$, such that
\begin{equation*}
\begin{split}
    \frac{1}{\epsilon_{a_{k}}^{2}}\Big(\big|\epsilon_{a_{k}}x+x_{a_{k}}|^{2}-1\Big)^{2}
    &=\epsilon_{a_{k}}^{2}\Big[|x|^{2}+\big|\frac{x_{a_{k}}}{\epsilon_{a_{k}}}\big|^{2}
    +2|x|\big|\frac{x_{a_{k}}}{\epsilon_{a_{k}}}\big|\cos\langle x,x_{a_{k}}\rangle
    -\frac{1}{\epsilon_{a_{k}}^{2}}\Big]^{2}\\
    &\geq\epsilon_{a_{k}}^{2}\Big(2|x|\big|\frac{x_{a_{k}}}{\epsilon_{a_{k}}}\big|\cos\langle x,x_{a_{k}}\rangle\Big)^{2}\\
    &\geq4|x_{a_{k}}|^{2}|x|^{2}\cos^{2}\big(\frac{\pi}{2}-\delta\big)\\
    &\geq C|x|^2\cos^{2}\big(\frac{\pi}{2}-\delta\big)\ \ \hbox{as}\ \ k\rightarrow\infty,
\end{split}
\end{equation*}
due to the fact that $\lim\limits_{k\rightarrow\infty}|x_{a_{k}}|= 1$. Taking $\delta=\frac{\pi}{20}$, the above estimate thus implies that
\begin{equation}\label{2.02}
\begin{split}
    &\lim\limits_{k\rightarrow\infty}\frac{1}{\epsilon_{a_{k}}^{2}}\int_{B_{R}(0)}\Big(|\epsilon_{a_{k}}x+x_{a_{k}}|^{2}-1\Big)^{2}|w_{a_{k}}|^{2}dx\\
    \geq& C \cos^{2}\frac{9\pi}{20}\int_{\Omega_{2}}|x|^{2}|w_{0}|^{2}dx
    \geq C(R),
\end{split}
\end{equation}
where $w_{0}>0$ is as in \eqref{conv7}.
The combination of \eqref{2.01} and \eqref{2.02} therefore yields that the proof of Lemma \ref{lem2.6} is complete.
\qed
\vskip 0.05truein

\noindent\textbf{Proof of Theorem \ref{th2.1}.}
By Lemma \ref{lem2.1}, it suffices to prove that there exists a constant $C>0$, independent of $a$, such that
\begin{equation}\label{low}
    I(a,\Omega_{a})\geq C(a^*-a)^{\frac{1}{2}-\beta}\ \ \hbox{as}\ \ a\nearrow a^*,
\end{equation}
where $\Omega_{a}=C_{0}(a^*-a)^{-\beta}$ holds for some $\beta\in[0,\frac{1}{2})$ and $C_{0}>0$.

On the contrary, assume that there exists a subsequence $\{a_{k}\}$, where $a_{k}\nearrow a^*$ as $k\to\infty$, such that
\begin{equation}\label{2:10}
    I(a_{k},\Omega_{a_{k}})=o\big((a^*-a_{k})^{\frac{1}{2}-\beta}\big)\ \ \hbox{as}\ \ k\to\infty.
\end{equation}
By the proof of Lemma \ref{lem2.3} (3), we obtain that there exists a convergent subsequence, still denoted by $\{a_{k}\}$, of $\{a_{k}\}$ such that $|w_{a_{k}}|\rightarrow \frac{Q}{\sqrt{a^*}}>0$ strongly in $L^{4}(\R^2)$ as $k\to\infty$. This implies that there exists a constant $M_{1}>0$, independent of $a_{k}$, such that
\begin{equation}\label{2:11}
    \int_{\R^2}|w_{a_{k}}|^{4}dx\geq M_{1}\ \ \hbox{as}\ \ k\rightarrow\infty.
\end{equation}
Moreover,  we derive from Lemma \ref{lem2.6} that there exists a constant $M_{2}>0$ such that
\begin{equation}\label{pot}
    \frac{\Omega_{a_{k}}^{2}}{8}\int_{\R^2}\big(|\epsilon_{a_{k}}x+x_{a_{k}}|^{2}-1\big)^{2}|w_{a_{k}}|^{2}dx
    \geq M_{2}\Omega_{a_{k}}^{2}\epsilon_{a_{k}}^{2}\ \ \hbox{as}\ \ k\rightarrow\infty.
\end{equation}
We then deduce from (\ref{2:10})--(\ref{pot}) that
\begin{equation}\label{low1}
\begin{aligned}
  &\quad o\big((a^*-a_{k})^{\frac{1}{2}-\beta}\big)\\
   &= I(a_{k},\Omega_{a_k{}})=F_{a_{k},\Omega_{a_{k}}}(u_{a_{k}})\\
    &\geq\int_{\R^2}\big|\nabla|u_{a_{k}}|\big|^{2}dx
    +\frac{\Omega_{a_{k}}^{2}}{8}\int_{\R^2}(|x|^{2}-1)^{2}|u_{a_{k}}|^{2}dx
    -\frac{a_{k}}{2}\int_{\R^2}|u_{a_{k}}|^{4}dx\\
    &=\frac{1}{\epsilon_{a_{k}}^{2}}\int_{\R^2}\big|\nabla|w_{a_{k}}|\big|^{2}dx
    +\frac{\Omega_{a_{k}}^{2}}{8}\int_{\R^2}(|\epsilon_{a_{k}}x+x_{a_{k}}|^{2}-1)^{2}|w_{a_{k}}|^{2}dx
    -\frac{a_{k}}{2\epsilon_{a_{k}}^{2}}\int_{\R^2}|w_{a_{k}}|^{4}dx\\
    &=\frac{1}{\epsilon_{a_{k}}^{2}}\Big(\int_{\R^2}\big|\nabla|w_{a_{k}}|\big|^{2}dx
    -\frac{a^*}{2}\int_{\R^2}|w_{a_{k}}|^{4}dx\Big)
    +\frac{a^*-a_{k}}{2\epsilon_{a_{k}}^{2}}\int_{\R^2}|w_{a_{k}}|^{4}dx\\
    &\quad+\frac{\Omega_{a_{k}}^{2}}{8}\int_{\R^2}\big(|\epsilon_{a_{k}}x+x_{a_{k}}|^{2}-1\big)^{2}|w_{a_{k}}|^{2}dx\\
    &\geq\frac{a^*-a_{k}}{2\epsilon_{a_{k}}^{2}}M_{1}+M_{2}\Omega_{a_{k}}^{2}\epsilon_{a_{k}}^{2}
    \geq C_{0}\sqrt{2M_{1}M_{2}}(a^*-a_{k})^{\frac{1}{2}-\beta}\ \ \hbox{as}\ \ k\rightarrow\infty,
\end{aligned}
\end{equation}
a contradiction, and hence (\ref{low}) holds true.
This completes the proof of Theorem \ref{th2.1}.\qed

Applying Theorem \ref{th2.1}, we also obtain from \eqref{low1} that
\begin{equation}\label{2.03}
    \epsilon_{a}=O\Big((a^*-a)^{\frac{1+2\beta}{4}}\Big)\ \ \hbox{and}\ \ |x_{a}|^2-1=O(\epsilon_{a})\ \ \hbox{as}\ \ a\nearrow a^*.
\end{equation}
Following (\ref{2.03}), we shall establish in next section the refined asymptotic estimates of complex-valued minimizers $u_{a}$ for $I(a,\Omega_{a})$ as $a\nearrow a^*$.

\section{Limiting Behavior of Minimizers as $a\nearrow a^*$}

Applying those results of previous sections, this section is devoted to the proof of Theorem \ref{th1}. We begin with the following $L^{\infty}-$uniform convergence of $w_{a_{k}}(x)$ as $k\to\infty$.

\vskip 0.05truein

\begin{prop}\label{pro}
Under the assumptions of Lemma \ref{lem2.3}, let $\{w_{a_k}\}$ be the subsequence obtained in Lemma \ref{lem2.3} (3).
Then we have
\begin{enumerate}
  \item
  There exists a large constant $R>0$ such that
  \begin{equation}\label{expo1}
    |w_{a_{k}}(x)|\leq Ce^{-\frac{2}{3}|x|}\ \ \hbox{in}\ \ \R^2\backslash B_{R}(0)\ \ \hbox{as}\ \ k\to\infty.
  \end{equation}
  \item The global maximum point $x_{a_{k}}$ of $|u_{a_{k}}|$ is unique as $k\to\infty$, and $w_{a_{k}}(x)$ satisfies
  \begin{equation}\label{wauni}
    w_{a_{k}}(x)\rightarrow\frac{Q(x)}{\sqrt{a^*}}\ \ \hbox{uniformly in}\ \ L^{\infty}(\R^2,\C)\ \ \hbox{as}\ \ k\to\infty.
  \end{equation}

  \item The following estimate holds:
  \begin{equation}\label{xz}
    \Omega_{a_{k}}\int_{\R^2}x^{\bot}\cdot\big(iw_{a_{k}}(x),\nabla w_{a_{k}}(x)\big)dx=o\big(\Omega_{a_{k}}^{2}\epsilon_{a_{k}}^{2}\big)\ \ \hbox{as}\ \ k\to\infty,
  \end{equation}
  where  $\epsilon_{a_{k}}>0$ is as in \eqref{ypx7}.
\end{enumerate}

\end{prop}

\noindent\textbf{Proof.}
1. Similar to the proof of \cite[Proposition 3.3]{GLY}, the exponential decay of $w_{a_{k}}$ in (\ref{expo1}) holds true and we omit the detailed proof for simplicity.

2. We first prove that
\begin{equation}\label{3.1}
    \lim\limits_{k\to\infty}w_{a_{k}}(x)\to\frac{Q(x)}{\sqrt{a^*}}\ \ \hbox{strongly in}\ \ H^1(\R^2,\C).
\end{equation}
By the definition of $w_{a_{k}}$ in (\ref{wa1}), we derive from \eqref{2.03} and \eqref{expo1} that
\begin{equation}\label{I1}
\begin{split}
    I(a_{k},\Omega_{a_{k}})&=F_{a_{k},\Omega_{a_{k}}}(u_{a_{k}})
    =F_{a_{k},\Omega_{a_{k}}}\Big(\frac{1}{\epsilon_{a_{k}}}e^{i(\Omega_{a_{k}} x\cdot x_{a_{k}}^{\bot}-\theta_{a_{k}})}w_{a_{k}}\big(\frac{x-x_{a_{k}}}{\epsilon_{a_{k}}}\big)\Big)\\
    &=\frac{1}{\epsilon_{a_{k}}^{2}}\int_{\R^2}|\nabla w_{a_{k}}|^{2}dx
    +\Omega_{a_{k}}^{2}\epsilon_{a_{k}}^{2}\int_{\R^2}|x|^{2}|w_{a_{k}}|^{2}dx\\
    &\quad-2\Omega_{a_{k}}\int_{\R^2}x^{\bot}\cdot(iw_{a_{k}},\nabla w_{a_{k}})dx
    -\frac{a_{k}}{2\epsilon_{a_{k}}^{2}}\int_{\R^2}|w_{a_{k}}|^{4}dx\\
    &\quad+\frac{\Omega_{a_{k}}^{2}}{8}\int_{\R^2}\big(|\epsilon_{a_{k}}x+x_{a_{k}}|^{2}-1\big)^{2}
    |w_{a_{k}}|^{2}dx\\
    &=\frac{1}{\epsilon_{a_{k}}^{2}}\Big(\int_{\R^2}|\nabla w_{a_{k}}|^{2}dx-\frac{a_{k}}{2}\int_{\R^2}|w_{a_{k}}|^{4}dx\Big)+o(\epsilon_{a_{k}}^{-2})\ \ \hbox{as}\ \ k\to\infty.
\end{split}
\end{equation}
Applying Theorem \ref{th2.1}, we deduce from \eqref{GN} and (\ref{I1}) that
\begin{equation*}
\begin{split}
    &\lim\limits_{k\to\infty}\int_{\R^2}\Big(|\nabla w_{a_{k}}|^{2}-\frac{a_{k}}{2}|w_{a_{k}}|^{4}\Big)dx\\
    =&\lim\limits_{k\to\infty}\int_{\R^2}\Big(|\nabla w_{a_{k}}|^{2}-\big|\nabla |w_{a_{k}}|\big|^{2}+\big|\nabla |w_{a_{k}}|\big|^{2}-\frac{a_{k}}{2}|w_{a_{k}}|^{4}\Big)dx=0,\\
\end{split}
\end{equation*}
which implies that
\begin{equation}\label{1.011}
    \int_{\R^2}|\nabla w_{a_{k}}|^{2}dx=\int_{\R^2}\big|\nabla |w_{a_{k}}|\big|^{2}dx+o(1)=\frac{a_{k}}{2}\int_{\R^2}|w_{a_{k}}|^{4}dx+o(1)
    \ \ \hbox{as}\ \ k\to\infty.
\end{equation}
By Lemma \ref{lem2.3}, we obtain from \eqref{1.011} that
\begin{equation}\label{wacon1}
    \lim\limits_{k\to\infty}w_{a_{k}}(x)=\frac{Q(x)}{\sqrt{a^*}}e^{i\sigma}\ \ \hbox{strongly in}\ \ H^{1}(\R^2,\C)
\end{equation}
for some $\sigma\in\R$.
We further have $\sigma=0$ in (\ref{wacon1}) in view of (\ref{orth4}), which implies that \eqref{3.1} is proved.

We next claim that $w_{a_{k}}(x)$ converges to $\frac{Q(x)}{\sqrt{a^*}}$ uniformly in $L^{\infty}(\R^2,\C)$ as $k\to\infty$.
Indeed, by the exponential decay of (\ref{expo}) and (\ref{expo1}), we only need to show the $L^{\infty}-$ uniform convergence of $w_{a_{k}}(x)$ on any compact domain of $\R^2$ as $k\to\infty$. Since $w_{a_{k}}(x)$ satisfies (\ref{waele}), denote
\begin{equation*}
    G_{a_{k}}(x):=-2i\Omega_{a_{k}}\epsilon_{a_{k}}^{2}(x^{\bot}\cdot\nabla w_{a_{k}})-\Big[\Omega_{a_{k}}^{2}\epsilon_{a_{k}}^{4}|x|^{2}
    +\frac{\Omega_{a_{k}}^{2}\epsilon_{a_{k}}^{2}}{8}\big(|\epsilon_{a_{k}}x+x_{a_{k}}|^{2}-1\big)^{2}
    -\mu_{a_{k}}\epsilon_{a_{k}}^{2}-a|w_{a_{k}}|^{2}\Big]w_{a_{k}},
\end{equation*}
so that
\begin{equation}\label{Ga1}
    -\Delta w_{a_{k}}(x)=G_{a_{k}}(x)\ \ \hbox{in}\ \ \R^2.
\end{equation}
Because $\{w_{a_{k}}(x)\}$ is bounded uniformly in $H^{1}(\R^2,\C)$ as $k\to\infty$, $\{G_{a_{k}}(x)\}$ is also bounded uniformly in $L^{2}_{loc}(\R^2,\C)$ as $k\to\infty$, where the estimates \eqref{2.03} and \eqref{wacon1} are also used. For any large $R>0$, we obtain from \cite[Theorem 8.8]{GT} that
\begin{equation}\label{8}
    \|w_{a_{k}}(x)\|_{H^{2}(B_{R})}\leq C\Big(\|w_{a_{k}}(x)\|_{H^{1}(B_{R+1})}+\|G_{a_{k}}(x)\|_{L^{2}(B_{R+1})}\Big),
\end{equation}
where $C>0$ is independent of $k>0$ and $R>0$. Therefore, $\{w_{a_{k}}(x)\}$ is also bounded uniformly in $H^{2}_{loc}(\R^2,\C)$ as $k\to\infty$. Since the embedding $H^{2}(B_{R})\hookrightarrow L^{\infty}(B_{R})$ is compact, cf. \cite[Theorem 7.26]{GT}, there exist a subsequence, still denoted by $\{w_{a_{k}}\}$, of $\{w_{a_{k}}\}$ and a function $\tilde{w}\in L^{\infty}(B_{R},\C)$ such that
\begin{equation*}
    w_{a_{k}}(x)\rightarrow \tilde{w}(x)\ \ \hbox{uniformly in}\ \ L^{\infty}(B_{R},\C)\ \ \hbox{as}\ \ k\rightarrow\infty,
\end{equation*}
where $a_{k}\nearrow a^*$ as $k\rightarrow\infty$. By (\ref{3.1}), we have $\tilde{w}(x)\equiv\frac{Q(x)}{\sqrt{a^*}}$.
Since $R>0$ is arbitrary, we get that
\begin{equation}\label{wauni1}
    w_{a_{k}}(x)\rightarrow\frac{Q(x)}{\sqrt{a^*}}\ \ \hbox{uniformly in}\ \ L^{\infty}_{loc}(\R^2,\C)\ \ \hbox{as}\ \ k\to\infty.
\end{equation}
Hence, (\ref{wauni}) holds in $L^{\infty}(\R^2,\C)$ as $k\to\infty$.

We next prove the uniqueness of $x_{a_{k}}$ as $k\to\infty$, where $x_{a_{k}}$ is a global maximum point of $|u_{a_{k}}|$. It follows from (\ref{8}) that $\{|\nabla w_{a_{k}}|\}$ is bounded uniformly in $L^{q}_{loc}(\R^2)$ as $k\to\infty$ for any $q\geq2$. The $L^{p}$ estimate \cite[Theorem 9.11]{GT} applied to (\ref{Ga1}) then implies that $\{w_{a_{k}}\}$ is bounded uniformly in $W^{2,q}_{loc}(\R^2)$ as $k\to\infty$. The standard Sobolev embedding theorem gives that $\{w_{a_{k}}\}$ is bounded uniformly in $C^{1,\alpha}_{loc}(\R^2)$. It then follows from the Schauder estimate \cite[Theorem 6.2]{GT} that $\{w_{a_{k}}\}$ is bounded uniformly in $C^{2,\alpha}_{loc}(\R^2)$ as $k\to\infty$. Hence, there exists a function $w_{1}\in C^{2}_{loc}(\R^2)$ such that
\begin{equation*}
    w_{a_{k}}\rightarrow w_{1}\ \ \hbox{in}\ \ C^{2}_{loc}(\R^2)\ \ \hbox{as}\ \ k\to\infty,
\end{equation*}
and we further deduce from (\ref{wauni}) that
\begin{equation*}
    |w_{a_{k}}(x)|\rightarrow w_{1}(x)\equiv\frac{Q(x)}{\sqrt{a^*}}\ \ \hbox{in}\ \ C^{2}_{loc}(\R^2)\ \ \hbox{as}\ \ k\to\infty.
\end{equation*}
Because the origin is the unique global maximum point of $Q(x)$, the above convergence shows that all global maximum points of $|w_{a_{k}}|$ must stay in a small ball $B_{\delta}(0)$ as $k\to\infty$ for some $\delta>0$. Since $Q''(0)<0$, we conclude that $Q''(r)<0$ for all $0<r<\delta$. It then follows from \cite[Lemma 4.2]{NT} that each $|w_{a_{k}}|$ has a unique global maximum point as $k\to\infty$, which is just the origin. This further proves the uniqueness of global maximum points $x_{a_{k}}$ as $k\to\infty$.

3. 
Denoting $w_{a_{k}}(x):=R_{a_{k}}(x)+iI_{a_{k}}(x)$, it then follows from (\ref{expo1}) and (\ref{wauni}) that there exists a constant $C>0$ such that
\begin{equation}\label{rot1}
\begin{split}
    \int_{\R^2}x^{\bot}\cdot(iw_{a_{k}},\nabla w_{a_{k}})dx
    &=\int_{\R^2}x^{\bot}\cdot(R_{a_{k}}\nabla I_{a_{k}}-I_{a_{k}}\nabla R_{a_{k}})dx\\
    &=2\int_{\R^2}x^{\bot}\cdot(R_{a_{k}}\nabla I_{a_{k}})dx\leq C\|\nabla I_{a_{k}}\|_{L^{2}(\R^2)}.
\end{split}
\end{equation}
We then derive from \eqref{I1} and (\ref{rot1}) that
\begin{equation}\label{ypx4}
\begin{split}
    \epsilon_{a_{k}}^{2}I(a_{k},\Omega_{a_{k}})&\geq\int_{\R^2}\big(|\nabla R_{a_{k}}|^{2}+|\nabla I_{a_{k}}|^{2}\big)dx\\
    &\quad-\frac{a^*}{2}\int_{\R^2}\big(R_{a_{k}}^4+I_{a_{k}}^{4}+2R_{a_{k}}^{2}I_{a_{k}}^{2}\big)dx
    -C\Omega_{a_{k}}\epsilon_{a_{k}}^{2}\|\nabla I_{a_{k}}\|_{L^{2}(\R^2)}.
\end{split}
\end{equation}
 Note from (\ref{wauni}) that $R_{a_{k}}\rightarrow\frac{Q}{\sqrt{a^*}}$ and $I_{a_{k}}\rightarrow 0$ uniformly in $\R^2$ as $k\to\infty$. We thus derive from \eqref{GN} that
\begin{equation}\label{est3}
\begin{split}
    &\quad\int_{\R^2}|\nabla R_{a_{k}}|^{2}dx-\frac{a^*}{2}\int_{\R^2}(R_{a_{k}}^{4}+I_{a_{k}}^{4}
    +2R_{a_{k}}^{2}I_{a_{k}}^{2})dx\\
    &\geq\int_{\R^2}|\nabla R_{a_{k}}|^{2}dx\Big(1-\int_{\R^2}R_{a_{k}}^{2}dx\Big)
    -\frac{a^*}{2}\int_{\R^2}(I_{a_{k}}^{2}+2R_{a_{k}}^{2})I_{a_{k}}^{2}dx\\
    &=\big(1+o(1)\big)\int_{\R^2}I_{a_{k}}^{2}dx-\big(1+o(1)\big)\int_{\R^2}Q^{2}I_{a_{k}}^{2}dx\ \ \hbox{as}\ \ k\to\infty.
\end{split}
\end{equation}
It then follows from (\ref{ypx4}) and (\ref{est3}) that
\begin{equation}\label{ypx5}
\begin{split}
    \epsilon_{a_{k}}^{2}I(a_{k},\Omega_{a_{k}})&\geq\int_{\R^2}|\nabla I_{a_{k}}|^{2}dx+\int_{\R^2}I_{a_{k}}^{2}dx\\
    &\quad-\int_{\R^2}Q^{2}I_{a_{k}}^{2}dx-o(1)\int_{\R^2}I_{a_{k}}^{2}dx
    -C\Omega_{a_{k}}\epsilon_{a_{k}}^{2}\|\nabla I_{a_{k}}\|_{L^{2}(\R^2)}\\
    &=(\mathcal{L}I_{a_{k}},I_{a_{k}})-o(1)\int_{\R^2}I_{a_{k}}^{2}dx-C\Omega_{a_{k}}\epsilon_{a_{k}}^{2}\|\nabla I_{a_{k}}\|_{L^{2}(\R^2)}\ \ \hbox{as}\ \ k\to\infty,
\end{split}
\end{equation}
where the operator $\mathcal{L}$ is defined by
\begin{equation}\label{L}
    \mathcal{L}:=-\Delta+1-Q^{2}\ \ \hbox{in}\ \ \R^{2}.
\end{equation}

We note from \cite[Corollary 11.9 and Theorem 11.8]{Lieb} that $(0,Q)$ is the first eigenpair of $\mathcal{L}$ and $ker(\mathcal{L})=\{Q\}$. Furthermore, since the essential spectrum of $\mathcal{L}$ satisfies $\sigma_{ess}(\mathcal{L})=[1,+\infty)$, we have $0\in\sigma_{d}(\mathcal{L})$, where $\sigma_{d}$ denotes the discrete spectrum of $\mathcal{L}$. We obtain from (\ref{orth4}) that
\begin{equation}\label{L1}
    (\mathcal{L}I_{a_{k}},I_{a_{k}})\geq\|I_{a_{k}}\|_{L^{2}(\R^2)}^{2}\ \ \hbox{as}\ \ k\to\infty.
\end{equation}
On the other hand, it follows from (\ref{L}) that
\begin{equation}\label{L2}
    (\mathcal{L}I_{a_{k}},I_{a_{k}})\geq\|\nabla I_{a_{k}}\|_{L^{2}(\R^2)}^{2}-\|Q\|^{2}_{L^{\infty}}\|I_{a_{k}}\|_{L^{2}(\R^2)}^{2}.
\end{equation}
Combining (\ref{L1}) with (\ref{L2}), there exists a constant $\bar{\rho}>0$ such that
\begin{equation}\label{L3}
    (\mathcal{L}I_{a_{k}},I_{a_{k}})\geq\bar{\rho}\|I_{a_{k}}\|^{2}_{H^{1}(\R^2)}\ \ \hbox{as}\ \ k\to\infty.
\end{equation}
Following Theorem \ref{th2.1}, we obtain from \eqref{2.03}, (\ref{ypx5}) and (\ref{L3}) that
\begin{equation*}
    C\Omega_{a_{k}}^2\epsilon^{4}_{a_{k}}\geq\epsilon_{a_{k}}^{2}I(a_{k},\Omega_{a_{k}})
    \geq\frac{\bar{\rho}}{2}\|I_{a_{k}}\|^{2}_{H^{1}(\R^2)}-C\Omega_{a_{k}}\epsilon_{a_{k}}^{2}\|\nabla I_{a_{k}}\|_{L^{2}(\R^2)}\ \ \hbox{as}\ \ k\to\infty,
\end{equation*}
which gives that
\begin{equation}\label{est7}
    \|I_{a_{k}}\|_{H^{1}(\R^2)}\leq C\Omega_{a_{k}}\epsilon_{a_{k}}^{2}\ \ \hbox{as}\ \ k\to\infty.
\end{equation}
Applying (\ref{expo1}) and (\ref{est7}), it yields from \eqref{2.03} that
\begin{equation*}
\begin{split}
    &\quad\Omega_{a_{k}}\int_{\R^2}x^{\bot}\cdot\big(iw_{a_{k}}(x),\nabla w_{a_{k}}(x)\big)dx\\
    &=2\Omega_{a_{k}}\int_{\R^2}x^{\bot}\cdot(R_{a_{k}}\nabla I_{a_{k}})dx=2\Omega_{a_{k}}\int_{\R^2}x^{\bot}\cdot\Big(\frac{Q(x)}{\sqrt{a^*}}\nabla I_{a_{k}}\Big)dx+o(\Omega_{a_{k}}^{2}\epsilon_{a_{k}}^{2})\\
    &=-2\Omega_{a_{k}}\int_{\R^2}x^{\bot}\cdot\Big(I_{a_{k}}\nabla\frac{Q(x)}{\sqrt{a^*}}\Big)dx
    +o(\Omega_{a_{k}}^{2}\epsilon_{a_{k}}^{2})\\
    &=o(\Omega_{a_{k}}^{2}\epsilon_{a_{k}}^{2})
    =o\big((a^*-a_{k})^{\frac{1}{2}-\beta}\big)\ \ \hbox{as}\ \ k\to\infty,
\end{split}
\end{equation*}
where the $L^{\infty}-$uniform convergence of $\{w_{a_{k}}\}$ in (\ref{wauni}) is also used. The proof of Proposition \ref{pro} is therefore complete.
\qed

\vskip 0.05truein
\noindent\textbf{Proof of Theorem \ref{th1}.}
By \eqref{wa1} and Proposition \ref{pro} (2), it remains to prove that
\begin{equation*}
    \epsilon_{a_{k}}:=\Big(\int_{\R^2}\big|\nabla |u_{a_{k}}|\big|^{2}dx\Big)^{-\frac{1}{2}}
    =\frac{(a^*-a_{k})^{\frac{1+2\beta}{4}}}{\sqrt{C_{0}}\lambda}+o\big((a^*-a_{k})^{\frac{1+2\beta}{4}}\big)>0\ \ \hbox{as}\ \ k\to\infty,
\end{equation*}
where $\lambda>0$ is defined as in (\ref{lmd}).

Actually, it follows from (\ref{I1}) that
\begin{equation}\label{Ia}
\begin{split}
    I(a_{k},\Omega_{a_{k}})=F_{a_{k},\Omega_{a_{k}}}(u_{a_{k}})&
    =\frac{1}{\epsilon_{a_{k}}^{2}}\Big[\int_{\R^2}|\nabla w_{a_{k}}|^{2}dx
    -\frac{a^*}{2}\int_{\R^2}|w_{a_{k}}|^{4}dx\Big]\\
    &\quad+\frac{a^*-a_{k}}{2\epsilon_{a_{k}}^{2}}\int_{\R^2}|w_{a_{k}}|^{4}dx
    +\Omega_{a_{k}}^{2}\epsilon_{a_{k}}^{2}\int_{\R^2}|x|^{2}|w_{a_{k}}|^{2}dx\\
    &\quad+\frac{\Omega_{a_{k}}^{2}}{8}\int_{\R^2}\big(|\epsilon_{a_{k}}x+x_{a_{k}}|^{2}-1\big)^{2}|w_{a_{k}}|^{2}dx\\
    &\quad-2\Omega_{a_{k}}\int_{\R^2}x^{\bot}\cdot(iw_{a_{k}},\nabla w_{a_{k}})dx.
\end{split}
\end{equation}
Since the term in the square bracket is non-negative, we drop it for a lower bound of $I(a_{k},\Omega_{a_{k}})$. Recall from (\ref{xz}) that
\begin{equation}\label{3:7}
    2\Omega_{a_{k}}\int_{\R^2}x^{\bot}\cdot(iw_{a_{k}},\nabla w_{a_{k}})dx=o\big(\Omega_{a_{k}}^{2}\epsilon_{a_{k}}^{2}\big)\ \ \hbox{as}\ \ k\to\infty.\ \
\end{equation}
Since $w_{a_{k}}\rightarrow\frac{Q}{\sqrt{a^*}}$ uniformly in $L^\infty(\R^2,\C)$ as $k\to\infty$, we get from \eqref{Q1} and the exponential decay (\ref{expo1}) that
\begin{equation}\label{3:8}
\begin{split}
    &\frac{a^*-a_{k}}{2\epsilon_{a_{k}}^{2}}\int_{\R^2}|w_{a_{k}}|^{4}dx
    =\big[1+o(1)\big]\frac{a^*-a_{k}}{2(a^*)^{2}\epsilon_{a_{k}}^{2}}\int_{\R^2}Q^{4}dx
    =\big[1+o(1)\big]\frac{a^*-a_{k}}{\epsilon_{a_{k}}^{2}a^*}\ \ \hbox{as}\ \ k\to\infty,\\
    &\Omega_{a_{k}}^{2}\epsilon_{a_{k}}^{2}\int_{\R^2}|x|^{2}|w_{a_{k}}|^{2}dx
    =\big[1+o(1)\big]\frac{\Omega_{a_{k}}^{2}\epsilon_{a_{k}}^{2}}{a^*}\int_{\R^2}|x|^{2}Q^{2}dx\ \ \hbox{as}\ \ k\to\infty.\\
\end{split}
\end{equation}
Note from Lemma \ref{lem2.3} (3) that $\lim\limits_{k\to\infty}x_{a_{k}}=x_{0}$ and $|x_{0}|=1$. We hence obtain from (\ref{2.03}) and (\ref{expo1}) that
\begin{equation}\label{3:6}
\begin{split}
    &\quad\frac{\Omega_{a_{k}}^{2}}{8}
    \int_{\R^2}\big(|\epsilon_{a_{k}}x+x_{a_{k}}|^{2}-1\big)^{2}|w_{a_{k}}|^{2}dx\\
    &=\frac{\Omega_{a_{k}}^{2}}{8}\int_{\R^2}\big[4\epsilon_{a_{k}}^2(x\cdot x_{a_{k}})^2+(|x_{a_{k}}|^2-1)^2+4\epsilon_{a_{k}}(|x_{a_{k}}|^2-1)(x\cdot x_{a_{k}})\big]|w_{a_{k}}|^{2}dx+o(\Omega_{a_{k}}^{2}\epsilon_{a_{k}}^{2})\\
    &=\frac{\Omega_{a_{k}}^{2}\epsilon_{a_{k}}^{2}}{2a^*}\int_{\R^2}(x\cdot x_{0})^{2}Q^{2}(x)dx+\frac{\Omega_{a_{k}}^{2}(|x_{a_{k}}|^2-1)^2}{8}
    +o\big(\Omega_{a_{k}}^{2}\epsilon_{a_{k}}^{2}\big)\\
    &\geq \frac{\Omega_{a_{k}}^{2}\epsilon_{a_{k}}^{2}}{2a^*}\int_{\R^2}(x\cdot x_{0})^{2}Q^{2}(x)dx+o\big(\Omega_{a_{k}}^{2}\epsilon_{a_{k}}^{2}\big)\\
    &=\frac{\Omega_{a_{k}}^{2}\epsilon_{a_{k}}^{2}}{4a^*}\int_{\R^2}|x|^{2}Q^{2}(x)dx
    +o\big(\Omega_{a_{k}}^{2}\epsilon_{a_{k}}^{2}\big)\ \ \hbox{as}\ \ k\to\infty,\\
\end{split}
\end{equation}
due to the fact that
\begin{equation}\label{3:13}
\begin{aligned}
  \int_{\R^2}(x\cdot x_{0})^{2}Q^{2}dx&=\int_{\R^2}(x_{1}x_{01}+x_{2}x_{02})^2Q^2(x)dx
  =\frac{1}{2}\int_{\R^2}|x|^2Q^2(x)dx,
\end{aligned}
\end{equation}
where $x_{0}=(x_{01},x_{02})\in\R^2$, and $Q(x)=Q(|x|)$ is also used.
We then deduce that
\begin{equation}\label{3:9}
\begin{split}
    &\quad\frac{a^*-a_{k}}{\epsilon_{a_{k}}^{2}a^*}
    +\frac{\Omega_{a_{k}}^{2}\epsilon_{a_{k}}^{2}}{a^*}\int_{\R^2}|x|^{2}Q^{2}dx
    +\frac{\Omega_{a_{k}}^{2}\epsilon_{a_{k}}^{2}}{4a^*}\int_{\R^2}|x|^{2}Q^{2}dx\\
    &=\frac{a^*-a_{k}}{\epsilon_{a_{k}}^{2}a^*}
    +\frac{5\Omega_{a_{k}}^{2}\epsilon_{a_{k}}^{2}}{4a^*}\int_{\R^2}|x|^{2}Q^{2}dx\\
    &=\frac{a^*-a_{k}}{\epsilon_{a_{k}}^{2}a^*}
    +\frac{\Omega_{a_{k}}^{2}\epsilon_{a_{k}}^{2}}{a^*}\lambda^{4}
    \geq\frac{2C_{0}(a^*-a_{k})^{\frac{1}{2}-\beta}}{a^*}\lambda^{2}\ \ \hbox{as}\ \ k\to\infty,
\end{split}
\end{equation}
by taking $\epsilon_{a_{k}}=\frac{(a^*-a_{k})^{\frac{1+2\beta}{4}}}{\sqrt{C_{0}}\lambda}>0$. This implies from above that
\begin{equation}\label{low2}
    \liminf_{k\to\infty}\frac{I(a_{k},\Omega_{a_{k}})}{(a^*-a_{k})^{\frac{1}{2}-\beta}}
    \geq\frac{2C_{0}\lambda^{2}}{a^*}.
\end{equation}

On the other hand, we consider the trial function
\begin{equation*}
    u_{\alpha}(x)=\frac{\alpha}{(a^*-a_{k})^{\frac{1+2\beta}{4}}\|Q\|_{2}}Q\Big(\frac{\alpha(x-y_{0})}{(a^*-a_{k})^{\frac{1+2\beta}{4}}}\Big)
    e^{i\Omega_{a_{k}}  x\cdot y_{0}^{\bot}},
\end{equation*}
where $y_{0}\in\R^2$ satisfies $|y_{0}|=1$ and $\alpha\in (0,\infty)$ is to be determined later.
By the definition of $u_{\alpha}$, we have
\begin{equation}\label{cal1}
\begin{aligned}
    &\quad\int_{\R^2}\big|(\nabla-i\Omega_{a_{k}} x^{\bot})u_{\alpha}\big|^{2}dx\\
    &=\frac{\alpha^{2}}{(a^*-a_{k})^{\frac{1+2\beta}{2}}a^*}\int_{\R^2}\Big|\nabla \Big(Q\big(\frac{\alpha(x-y_{0})}{(a^*-a_{k})^{\frac{1+2\beta}{4}}}\big)e^{i\Omega_{a_{k}} x\cdot y_{0}^{\bot}}\Big)\\
    &\qquad\qquad\qquad\qquad\qquad
    -i\Omega_{a_{k}} x^{\bot}Q\Big(\frac{\alpha(x-y_{0})}{(a^*-a_{k})^{\frac{1+2\beta}{4}}}\Big)e^{i\Omega_{a_{k}} x\cdot y_{0}^{\bot} }\Big|^2dx\\
    &=\frac{\alpha^{2}}{(a^*-a_{k})^{\frac{1+2\beta}{2}}a^*}\int_{\R^2}\Big|\nabla Q-i\Omega_{a_{k}}\frac{(a^*-a_{k})^{\frac{1+2\beta}{2}}}{\alpha^{2}} x^{\bot}Q\Big|^2dx\\
    &=\frac{\alpha^{2}}{(a^*-a_{k})^{\frac{1+2\beta}{2}}}
    +\frac{(a^*-a_{k})^{\frac{1+2\beta}{2}}\Omega_{a_{k}}^{2}}{\alpha^{2}a^*}\int_{\R^2}|x|^{2}Q^2dx,\\
\end{aligned}
\end{equation}
and
\begin{equation}\label{cal2}
\begin{split}
    \quad\frac{a_{k}}{2}\int_{\R^2}|u_{\alpha}|^{4}dx
    &=\frac{\alpha^2a_{k}}{2(a^*)^{2}(a^*-a_{k})^{\frac{1+2\beta}{2}}}\int_{\R^2}Q^{4}dx
    =\frac{\alpha^2a_{k}}{a^*(a^*-a_{k})^{\frac{1+2\beta}{2}}}.
\end{split}
\end{equation}
We also deduce that
\begin{equation}\label{cal3}
\begin{split}
    &\quad\frac{\Omega_{a_{k}}^{2}}{8}\int_{\R^2}\big(|x|^{2}-1\big)^{2}|u_{\alpha}|^{2}dx\\
    &=\frac{\Omega_{a_{k}}^{2}}{8a^*}\int_{\R^2}
    \Big(\Big|\frac{(a^*-a_{k})^{\frac{1+2\beta}{4}}}{\alpha}x+y_{0}\Big|^{2}-1\Big)^{2}Q^2dx\\
    &=\frac{\Omega_{a_{k}}^{2}(a^*-a_{k})^{\frac{1+2\beta}{2}}}{4a^*\alpha^{2}}
    \int_{\R^2}|x|^{2}Q^{2}dx+o\big((a^*-a_{k})^{\frac{1}{2}-\beta}\big)\ \ \hbox{as}\ \ k\to\infty,
\end{split}
\end{equation}
where \eqref{3:13} is also used in the last equality.
It then follows from (\ref{cal1})--(\ref{cal3}) that
\begin{equation}\label{3:12}
\begin{aligned}
    I(a_{k},\Omega_{a_{k}})&\leq F_{a_{k},\Omega_{a_{k}}}(u_{\alpha})\\
    &=\int_{\R^2}\big|(\nabla-i\Omega_{a_{k}} x^{\bot})u_{\alpha}\big|^{2}dx
    -\frac{a_{k}}{2}\int_{\R^2}|u_{\alpha}|^{4}dx
    +\frac{\Omega_{a_{k}}^{2}}{8}\int_{\R^2}\big(|x|^{2}-1\big)^{2}|u_{\alpha}|^{2}dx\\
    &=\frac{\alpha^{2}(a^*-a_{k})^{\frac{1}{2}-\beta}}{a^*}
    +\frac{5\Omega_{a_{k}}^{2}(a^*-a_{k})^{\frac{1+2\beta}{2}}}{4a^*\alpha^{2}}
    \int_{\R^2}|x|^{2}Q^{2}dx+o\big((a^*-a_{k})^{\frac{1}{2}-\beta}\big)\\
    &=\frac{\alpha^{2}(a^*-a_{k})^{\frac{1}{2}-\beta}}{a^*}
    +\frac{\Omega_{a_{k}}^{2}(a^*-a_{k})^{\frac{1+2\beta}{2}}}{a^*\alpha^{2}}\lambda^{4}+o\big((a^*-a_{k})^{\frac{1}{2}-\beta}\big)\\
    &\leq\frac{2C_{0}(a^*-a_{k})^{\frac{1}{2}-\beta}}{a^*}\lambda^{2}+o\big((a^*-a_{k})^{\frac{1}{2}-\beta}\big)\ \ \hbox{as}\ \ k\to\infty,
\end{aligned}
\end{equation}
by taking $\alpha=\sqrt{C_{0}}\lambda>0$, which implies that
\begin{equation}\label{hest1}
    \limsup_{k\to\infty}\frac{I(a_{k},\Omega_{a_{k}})}{(a^*-a_{k})^{\frac{1}{2}-\beta}}
    \leq\frac{2C_{0}\lambda^{2}}{a^*}.
\end{equation}

We now conclude from \eqref{3:9}, (\ref{low2}) and (\ref{hest1}) that
\begin{equation}\label{conver}
    \lim\limits_{k\to\infty}\frac{I(a_{k},\Omega_{a_{k}})}{(a^*-a_{k})^{\frac{1}{2}-\beta}}
    =\frac{2C_{0}\lambda^{2}}{a^*},
\end{equation}
and
\begin{equation*}
    \epsilon_{a_{k}}=\Big(\int_{\R^2}\big|\nabla |u_{a_{k}}|\big|^{2}\Big)^{-\frac{1}{2}}
    =\frac{(a^*-a_{k})^{\frac{1+2\beta}{4}}}{\sqrt{C_{0}}\lambda}+o\big((a^*-a_{k})^{\frac{1+2\beta}{4}}\big)>0\ \ \hbox{as}\ \ k\to\infty.
\end{equation*}
Following above estimates, we finally claim that
\begin{equation}\label{3:10}
    \lim\limits_{k\to\infty}\frac{|x_{a_{k}}|^{2}-1}{(a^*-a_{k})^{\frac{1+2\beta}{4}}}=0.
\end{equation}
On the contrary, assume that there exist a constant $C>0$ and a subsequence of $\{a_{k}\}$, still denoted by $\{a_{k}\}$, such that
\begin{equation*}
  \frac{(|x_{a_{k}}|^2-1)^2}{(a^*-a_{k})^{\frac{1+2\beta}{2}}}\geq C\ \ \hbox{as}\ \ k\to\infty.
\end{equation*}
It then follows from \eqref{Ia}--\eqref{3:6}, \eqref{3:9} and \eqref{hest1} that
\begin{equation*}
\begin{aligned}
  \frac{2C_{0}\lambda^2}{a^*}&\geq \lim\limits_{k\to\infty}\frac{I(a_{k},\Omega_{a_{k}})}{(a^*-a_{k})^{\frac{1}{2}-\beta}}\\
  &\geq\lim\limits_{k\to\infty}\frac{1}{(a^*-a_{k})^{\frac{1}{2}-
  \beta}}\Big[\frac{a^*-a_{k}}{\epsilon_{a_{k}}^2a^*}+\frac{\Omega_{a_{k}}^2\epsilon_{a_{k}}^2}{a^*}\int_{\R^2}|x|^2Q^2dx
  +\frac{\Omega_{a_{k}}^{2}\epsilon_{a_{k}}^{2}}{4a^*}\int_{\R^2}|x|^{2}Q^{2}dx\\
  &\quad+\frac{\Omega_{a_{k}}^{2}(|x_{a_{k}}|^2-1)^2}{8}\Big]
  \geq\frac{2C_{0}\lambda^2}{a^*}+\frac{CC_{0}^2}{8}> \frac{2C_{0}\lambda^2}{a^*},
\end{aligned}
\end{equation*}
a contradiction, and  the claim \eqref{3:10} hence holds true.
This completes the proof of Theorem \ref{th1}.
\qed

\section{Refined Estimates of Minimizers as $a\nearrow a^*$}
Under the assumption that $\Omega_{a}>0$ satisfies \eqref{0:1}, i.e., $\Omega_{a}=C_{0}(a^*-a)^{-\beta}$ holds for some $\beta\in[0,\frac{1}{2})$ and $C_{0}>0$, the purpose of this section is to derive the refined estimates of complex-valued minimizers $u_{a}$ for $I(a,\Omega_{a})$ as $a\nearrow a^*$. We start with some assumptions and notations.

Assume $u_{a_{k}}$ is a minimizer of $I(a_{k},\Omega_{a_{k}})$ and $x_{a_{k}}$ is a global maximum point of $|u_{a_{k}}|$. It then follows from Theorem \ref{th1} that there exists a subsequence, still denoted by $\{(u_{a_{k}}, x_{a_{k}})\}$,  of $\{(u_{a_{k}}, x_{a_{k}})\}$
such that
\begin{equation}\label{4:22}
  \lim\limits_{k\to\infty}x_{a_{k}}=x_{0}
\end{equation}
 for some point $x_{0}\in\R^2$ satisfying $|x_{0}|=1$. For this convergent subsequence, we define
\begin{equation}\label{4.1}
    \eps_{a_{k}}:=\frac{(a^*-a_{k})^{\frac{1+2\beta}{4}}}{\sqrt{C_{0}}\lambda}>0,\ \ \lambda:=\Big[\frac{5}{4}\int_{\R^2}|x|^{2}Q^{2}dx\Big]^{\frac{1}{4}}>0,
\end{equation}
and
\begin{equation*}
    v_{a_{k}}(x):=\eps_{a_{k}}\sqrt{a^*}u_{a_{k}}(\eps_{a_{k}}x+x_{a_{k}})e^{-i(\eps_{a_{k}}\Omega_{a_{k}} x\cdot x_{a_{k}}^{\bot}-\theta_{a_{k}})}=\tilde{R}_{a_{k}}(x)+i\tilde{I}_{a_{k}}(x),
\end{equation*}
where $\tilde{R}_{a_{k}}(x)$ and $\tilde{I}_{a_{k}}(x)$ denote the real and imaginary parts of $v_{a_{k}}(x)$, respectively, and the constant phase $\theta_{a_{k}}\in[0,2\pi)$ is chosen such that
\begin{equation}\label{orth6}
    \|v_{a_{k}}-Q\|_{L^{2}(\R^2)}=\min\limits_{\theta\in [0,2\pi)}\|e^{i\theta}\tilde{v}_{a_{k}}-Q\|_{L^{2}(\R^2)},
\end{equation}
where $\tilde{v}_{a_{k}}=\eps_{a_{k}}\sqrt{a^*}u_{a_{k}}(\eps_{a_{k}}x+x_{a_{k}})e^{-i(\eps_{a_{k}}\Omega_{a_{k}} x\cdot x_{a_{k}}^{\bot})}$. Note that (\ref{orth6}) essentially gives the orthogonality condition on $\tilde{I}_{a_{k}}(x)$:
\begin{equation*}
    \int_{\R^2}Q(x)\tilde{I}_{a_{k}}(x)dx=0.
\end{equation*}
Using the above notations, the proof of Theorem \ref{th1} yields that
\begin{equation*}
    v_{a_{k}}(x)\rightarrow Q(x)\ \ \hbox{uniformly in}\ \ L^{\infty}(\R^2,\C)\ \ \hbox{as}\ \ k\to\infty,
\end{equation*}
and the unique global maximum point $x_{a_{k}}\in\R^2$ of $|u_{a_{k}}|$ satisfies
\begin{equation*}
    \lim\limits_{k\to\infty}\frac{|x_{a_{k}}|^{2}-1}{\eps_{a_{k}}}=0.
\end{equation*}
Following (\ref{ELE1}), we deduce that $v_{a_{k}}$ satisfies the following elliptic equation
\begin{equation}\label{4.6}
\begin{split}
    &\quad-\Delta v_{a_{k}}+2i\Omega_{a_{k}}\eps_{a_{k}}^{2}x^{\bot}\cdot\nabla v_{a_{k}}+\Big[\Omega_{a_{k}}^{2}\eps_{a_{k}}^{4}|x|^{2}
    +\frac{\Omega_{a_{k}}^{2}\eps_{a_{k}}^{2}}{8}(|\eps_{a_{k}}x+x_{a_{k}}|^{2}-1)^2\Big]v_{a_{k}}\\
    &=\mu_{a_{k}}\eps_{a_{k}}^{2}v_{a_{k}}+\frac{a_{k}}{a^*}|v_{a_{k}}|^{2}v_{a_{k}}\ \ \hbox{in}\ \ \R^2,
\end{split}
\end{equation}
and
\begin{equation}\label{4.7}
    \lim\limits_{k\to\infty}\mu_{a_{k}}\eps_{a_{k}}^{2}=-1.
\end{equation}
\vskip 0.05truein
To obtain the refined profile of $u_{a_{k}}$ as $k\to\infty$, we now define
\begin{equation}\label{6.1}
    \bar{\varepsilon}_{a_{k}}=\sqrt{-\frac{1}{\mu_{a_{k}}}}>0\ \ \hbox{as}\ \ k\to\infty,
\end{equation}
where the Lagrange multiplier $\mu_{a_{k}}$ is as in (\ref{4.6}). It then follows from (\ref{4.1}) and (\ref{4.7}) that $\bar{\varepsilon}_{a_{k}}\rightarrow 0$ as $k\to\infty$. We also denote
\begin{equation}\label{6.2}
\begin{aligned}
    \nu_{a_{k}}(x):&=\bar{\varepsilon}_{a_{k}}\sqrt{a_{k}}u_{a_{k}}(\bar{\varepsilon}_{a_{k}}x
    +x_{a_{k}})e^{-i(\bar{\varepsilon}_{a_{k}}\Omega_{a_{k}} x\cdot x_{a_{k}}^{\bot}-\rho_{a_{k}})}\\
    &=\tilde{\mathrm{R}}_{a_{k}}(x)+i\mathrm{I}_{a_{k}}(x)=\big[\mathrm{R}_{a_{k}}(x)+Q(x)\big]
    +i\mathrm{I}_{a_{k}}(x),
\end{aligned}
\end{equation}
where $x_{a_{k}}$ is the unique global maximum point of $|u_{a_{k}}|$ as $k\to\infty$, $\tilde{\mathrm{R}}_{a_{k}}(x)$ and $\mathrm{I}_{a_{k}}(x)$ are the real and imaginary parts of $\nu_{a_{k}}(x)$, respectively. The constant phase $\rho_{a_{k}}\in\R$ satisfies
\begin{equation}\label{6.3}
    \|\nu_{a_{k}}-Q(x)\|_{L^2(\R^2)}=\min\limits_{\rho\in[0,2\pi)}\|e^{i\rho}\tilde{\nu}_{a_{k}}-Q\|_{L^2(\R^2)},
\end{equation}
where $\tilde{\nu}_{a_{k}}:=\bar{\varepsilon}_{a_{k}}\sqrt{a_{k}}u_{a_{k}}(\bar{\varepsilon}_{a_{k}}x+x_{a_{k}})e^{-i(\bar{\varepsilon}_{a_{k}}\Omega_{a_{k}} x\cdot x_{a_{k}}^{\bot})}$. Note from (\ref{6.3}) that the term $\mathrm{I}_{a_{k}}$ satisfies
\begin{equation}\label{6.4}
    \int_{\R^2}Q(x)\mathrm{I}_{a_{k}}(x)dx=0.
\end{equation}
Similar to Theorem \ref{th1}, one can derive that
\begin{equation}\label{6.5}
    \bar{\varepsilon}_{a_{k}}=\varepsilon_{a_{k}}\big[1+o(1)\big]\ \ \hbox{as}\ \ k\to\infty,\ \ \mathrm{R}_{a_{k}}(x)\rightarrow 0\ \ \hbox{and}\ \ \mathrm{I}_{a_{k}}(x)\rightarrow 0\ \ \hbox{uniformly in}\ \ L^{\infty}(\R^2,\R)\ \ \hbox{as}\ \ k\to\infty,
\end{equation}
where $\eps_{a_{k}}>0$ and $\bar{\eps}_{a_{k}}>0$ are as in \eqref{4.1} and (\ref{6.1}), respectively.

\subsection{Refined estimates of $\nu_{a_{k}}$ as $k\to\infty$}
 Following \eqref{6.3}--\eqref{6.5}, in this subsection we derive the refined estimates of $\nu_{a_{k}}$ defined by (\ref{6.2}) as $k\to\infty$. Denote the operator
\begin{equation}\label{4.8}
    \mathcal{L}:=-\Delta+1-Q^{2}\ \ \hbox{in}\ \ \R^2.
\end{equation}
By \cite[Theorem 11.8]{Lieb} and \cite[Corollary 11.9]{Lieb}, we have
\begin{equation}\label{4.9}
    ker\mathcal{L}=\{Q\}\ \ \hbox{and}\ \ \big<\mathcal{L}v,v\big>\geq 0\ \ \hbox{for all}\ \ v\in L^{2}(\R^2).
\end{equation}
We also define the linearized operator $\tilde{\mathcal{L}}$ by
\begin{equation}\label{4.10}
    \tilde{\mathcal{L}}:=-\Delta+1-3Q^{2}\ \ \hbox{in}\ \ \R^2,
\end{equation}
and it then follows from \cite{K,NT} that
\begin{equation}\label{4.11}
    ker\tilde{\mathcal{L}}=\Big\{\frac{\partial Q}{\partial x_{1}},\frac{\partial Q}{\partial x_{2}}\Big\}.
\end{equation}


Recall from (\ref{ELE1}) and \eqref{6.2} that $\nu_{a_{k}}$ satisfies
\begin{equation}\label{6.7}
\begin{split}
    &\quad-\Delta \nu_{a_{k}}+2i\Omega_{a_{k}}\bar{\varepsilon}_{a_{k}}^{2}(x^{\bot}\cdot\nabla \nu_{a_{k}})+\Big[\Omega_{a_{k}}^{2}\bar{\varepsilon}_{a_{k}}^{4}|x|^{2}
    +\frac{\Omega_{a_{k}}^{2}\bar{\varepsilon}_{a_{k}}^{2}}{8}\big(|\bar{\varepsilon}_{a_{k}}x+x_{a_{k}}|^{2}-1\big)^2\Big]\nu_{a_{k}}\\
    &=-\nu_{a_{k}}+|\nu_{a_{k}}|^{2}\nu_{a_{k}}\ \ \hbox{in}\ \ \R^2,
\end{split}
\end{equation}
where $\bar{\varepsilon}_{a_{k}}>0$ is as in \eqref{6.1}.
Denote the linear operator
\begin{equation}\label{6.8}
    \mathrm{L}_{a_{k}}:=-\Delta+\Big[\Omega_{a_{k}}^{2}\bar{\varepsilon}_{a_{k}}^{4}|x|^{2}
    +\frac{\Omega_{a_{k}}^{2}\bar{\varepsilon}_{a_{k}}^{2}}{8}
    \big(|\bar{\varepsilon}_{a_{k}}x+x_{a_{k}}|^{2}-1\big)^{2}+1-|\nu_{a_{k}}|^{2}\Big]\ \ \hbox{in}\ \ \R^2.
\end{equation}
It then follows from (\ref{6.2}), (\ref{6.4}) and (\ref{6.7}) that $\mathrm{I}_{a_{k}}$ satisfies
\begin{equation}\label{6.9}
    \mathrm{L}_{a_{k}}\mathrm{I}_{a_{k}}
    =-2\Omega_{a_{k}}\bar{\varepsilon}_{a_{k}}^{2}(x^{\bot}\cdot\nabla\mathrm{R}_{a_{k}})\ \ \hbox{in}\ \ \R^2, \ \ \int_{\R^2}Q\mathrm{I}_{a_{k}}dx=0,
\end{equation}
and $\mathrm{R}_{a_{k}}$ satisfies
\begin{equation}\label{6.10}
    \tilde{\mathrm{L}}_{a_{k}}\mathrm{R}_{a_{k}}:=\big[\mathrm{L}_{a_{k}}-Q^2
    -Q\tilde{\mathrm{R}}_{a_{k}}\big]\mathrm{R}_{a_{k}}=\mathrm{F}_{a_{k}}(x)\ \ \hbox{in}\ \ \R^2,
\end{equation}
where $\mathrm{F}_{a_{k}}(x)$ is defined by
\begin{equation}\label{6.11}
\begin{split}
    \mathrm{F}_{a_{k}}(x):
    &=-\Big[\Omega_{a_{k}}^{2}\bar{\varepsilon}_{a_{k}}^{4}|x|^2
    +\frac{\Omega_{a_{k}}^{2}\bar{\varepsilon}_{a_{k}}^{2}}{8}\big(|\bar{\varepsilon}_{a_{k}}x+x_{a_{k}}|^{2}-1\big)^2\Big]Q
    +2\Omega_{a_{k}}\bar{\varepsilon}_{a_{k}}^{2}(x^{\bot}\cdot\nabla \mathrm{I}_{a_{k}})+\mathrm{I}_{a_{k}}^{2}Q\\
    &=-\Omega_{a_{k}}^{2}\bar{\varepsilon}_{a_{k}}^{4}\Big\{|x|^{2}
    +\frac{1}{8}\Big[\bar{\varepsilon}_{a_{k}}^{2}|x|^{4}+4(x\cdot x_{a_{k}})^{2}+\Big(\frac{|x_{a_{k}}|^{2}-1}{\bar{\varepsilon}_{a_{k}}}\Big)^{2}
    +4\bar{\varepsilon}_{a_{k}}|x|^{2}(x\cdot x_{a_{k}})\\
    &\quad+2(|x_{a_{k}}|^{2}-1)|x|^{2}+4\frac{|x_{a_{k}}|^{2}-1}{\bar{\varepsilon}_{a_{k}}}(x\cdot x_{a_{k}})\Big]\Big\}Q
    +2\Omega_{a_{k}}\bar{\varepsilon}_{a_{k}}^{2}(x^{\bot}\cdot\nabla \mathrm{I}_{a_{k}})+\mathrm{I}_{a_{k}}^{2}Q.\\
\end{split}
\end{equation}

Similar to \cite[Proposition 2.2]{GLP1}, we deduce from (\ref{6.7}) that there exists a constant $C>0$, independent of $k>0$, such that
\begin{equation}\label{4.44}
    |\nabla \nu_{a_{k}}|\leq Ce^{-\frac{1}{2}|x|}\ \ \hbox{uniformly in}\ \ \R^2\ \ \hbox{as}\ \ k\to\infty.
\end{equation}
Since $(\nabla |\nu_{a_{k}}|)(0)\equiv 0$ holds for all $k>0$, we obtain from (\ref{6.2}) and (\ref{4.44}) that
\begin{equation}\label{4.14}
    \nabla \mathrm{R}_{a_{k}}(0)=-\frac{\mathrm{I}_{a_{k}}(0)\nabla \mathrm{I}_{a_{k}}(0)}{Q(0)+\mathrm{R}_{a_{k}}(0)}\rightarrow 0\ \ \hbox{as}\ \ k\to\infty.
\end{equation}
Using the same argument of \cite[Lemma 4.3]{GLY}, we then obtain from (\ref{6.9}), (\ref{4.44}) and (\ref{4.14}) that there exists a constant $C>0$, independent of $k>0$, such that
\begin{equation}\label{4.45}
    |\nabla \mathrm{I}_{a_{k}}(x)|,\ \ |\mathrm{I}_{a_{k}}(x)|\leq C\Omega_{a_{k}}\bar{\eps}_{a_{k}}^{2}e^{-\frac{1}{8}|x|}\ \ \hbox{uniformly in}\ \ \R^2\ \ \hbox{as}\ \ k\to\infty.
\end{equation}
We next discuss the refined estimates of $\mathrm{R}_{a_{k}}$ and $\mathrm{I}_{a_{k}}$.

\begin{lem}\label{lem4.4}
Let $u_{a_{k}}$ be a complex-valued minimizer of $I(a_{k},\Omega_{a_{k}})$, where $\Omega_{a_{k}}>0$ satisfies \eqref{0:1} for some $\beta\in [0,\frac{1}{2})$ and $C_{0}>0$. Suppose that  $x_{a_{k}}\in\R^2$ is the unique maximum point of $|u_{a_{k}}|$ as $k\to\infty$. Then we have
\begin{enumerate}
  \item There exists a constant $C>0$, which is independent of $k>0$, such that the imaginary part $\mathrm{I}_{a_{k}}$ of (\ref{6.2}) satisfies
      \begin{equation}\label{4.47}
        |\nabla \mathrm{I}_{a_{k}}(x)|,\ \ |\mathrm{I}_{a_{k}}(x)|\leq C\Omega_{a_{k}}^3\bar{\eps}_{a_{k}}^{6}e^{-\frac{1}{14} |x|}\ \ \hbox{uniformly in}\ \ \R^2\ \ \hbox{as}\ \ k\to\infty.
      \end{equation}
  \item The real part $\mathrm{R}_{a_{k}}$ of (\ref{6.2}) satisfies
  \begin{equation}\label{4.48}
    \mathrm{R}_{a_{k}}(x)=\Omega_{a_{k}}^{2}\bar{\eps}_{a_{k}}^{4}\Psi_{1}(x)
    +o(\Omega_{a_{k}}^{2}\bar{\eps}_{a_{k}}^{4})\ \ \hbox{in}\ \ \R^2\ \ \hbox{as}\ \ k\to\infty,
  \end{equation}
  where $\Psi_{1}(x)\in C^{2}(\R^2)\cap L^{\infty}(\R^2)$ solves uniquely the following problem
  \begin{equation}\label{4.49}
    \nabla\Psi_{1}(0)=0,\ \ \tilde{\mathcal{L}}\Psi_{1}(x)=-\Big[|x|^2+\frac{1}{2}(x\cdot x_{0})^{2}\Big]Q\ \ \hbox{in}\ \ \R^2,
  \end{equation}
  where $x_{0}$ is as in \eqref{4:22}.
\end{enumerate}
\end{lem}

\noindent\textbf{Proof.}
\emph{Step 1.} Denote
\begin{equation}\label{4.50}
    \mathcal{R}_{a_{k}}(x):=\mathrm{R}_{a_{k}}(x)-\Omega_{a_{k}}^{2}\bar{\eps}_{a_{k}}^{4}\Psi_{1}(x),
\end{equation}
where $\Psi_{1}(x)\in C^{2}(\R^2)\cap L^{\infty}(\R^2)$ is a solution of (\ref{4.49}). The uniqueness of $\Psi_{1}(x)$ follows from (\ref{4.11}) and the fact that $\nabla\Psi_{1}(0)=0$. Moreover, by the comparison principle, we derive from (\ref{expo}) and (\ref{4.49}) that
\begin{equation}\label{4:1}
    |\Psi_{1}(x)|\leq Ce^{-\delta |x|}\ \ \hbox{in}\ \ \R^2,\ \ \hbox{where}\ \ \frac{4}{5}<\delta<1.
\end{equation}
It follows from (\ref{4.14}), (\ref{4.49}) and \eqref{4.50} that $\mathcal{R}_{a_{k}}$ satisfies
\begin{equation}\label{4:2}
    \nabla\mathcal{R}_{a_{k}}(0)\rightarrow 0\ \ \hbox{as}\ \ k\to\infty.
\end{equation}
Following (\ref{6.10}) and (\ref{4.50}), we have
\begin{equation}\label{4:3}
    \tilde{\mathrm{L}}_{a_{k}}\mathcal{R}_{a_{k}}=\mathcal{F}_{a_{k}}(x)-(\tilde{\mathrm{L}}_{a_{k}}
    -\tilde{\mathcal{L}})\Omega_{a_{k}}^{2}\bar{\eps}_{a_{k}}^{4}\Psi_{1}(x)
    :=N_{a_{k}}(x),
\end{equation}
where the operator $\tilde{\mathrm{L}}_{a_{k}}$ is as in \eqref{6.10}, $\mathcal{F}_{a_{k}}(x)=\mathrm{F}_{a_{k}}(x)
-\Omega_{a_{k}}^{2}\bar{\eps}_{a_{k}}^{4}\Big\{-\big[|x|^2+\frac{1}{2}(x\cdot x_{0})^{2}\big]Q\Big\}$, and $F_{a_{k}}(x)$ is as in \eqref{6.10}.
It yields from (\ref{expo}), (\ref{6.11}) (\ref{4.45}) and (\ref{4:1}) that
\begin{equation}\label{4:4}
    \frac{\mathcal{F}_{a_{k}}(x)}{\Omega_{a_{k}}^{2}\bar{\eps}_{a_{k}}^{4}}\leq Ce^{-\frac{1}{10}|x|}\ \ \hbox{uniformly in}\ \ \R^2\ \ \hbox{as}\ \ k\to\infty,
\end{equation}
and
\begin{equation}\label{4:5}
    \frac{\Big|(\tilde{\mathrm{L}}_{a_{k}}
    -\tilde{\mathcal{L}})\Omega_{a_{k}}^{2}\bar{\eps}_{a_{k}}^{4}\Psi_{1}(x)
    \Big|}{\Omega_{a_{k}}^{2}\bar{\eps}_{a_{k}}^{4}}\leq C\delta_{a_{k}}e^{-\frac{1}{10}|x|}\ \ \hbox{uniformly in}\ \ \R^2\ \ \hbox{as}\ \ k\to\infty,
\end{equation}
where $\delta_{a_{k}}>0$ satisfies $\delta_{a_{k}}=o(1)$ as $k\to\infty$. We then derive from (\ref{4:4}) and (\ref{4:5}) that
\begin{equation}\label{4:7}
    \frac{|N_{a_{k}}(x)|}{\Omega_{a_{k}}^{2}\bar{\eps}_{a_{k}}^{4}}\leq Ce^{-\frac{1}{10}|x|}\ \ \hbox{uniformly in}\ \ \R^2\ \ \hbox{as}\ \ k\to\infty.
\end{equation}

We now claim that there exists a constant $C>0$, independent of $k>0$, such that
\begin{equation}\label{4:6}
    |\mathcal{R}_{a_{k}}(x)|\leq C\Omega_{a_{k}}^{2}\bar{\eps}_{a_{k}}^{4}\ \ \hbox{uniformly in}\ \ \R^2\ \ \hbox{as}\ \ k\to\infty.
\end{equation}
Instead, assume that the claim (\ref{4:6}) is false, i.e., $\lim\limits_{k\to\infty}\frac{\|\mathcal{R}_{a_{k}}\|_{L^{\infty}(\R^2)}}{\Omega_{a_{k}}^{2}\bar{\eps}_{a_{k}}^{4}}=\infty$.
Denote $\hat{\mathcal{R}}_{a_{k}}=\frac{\mathcal{R}_{a_{k}}}{\|\mathcal{R}_{a_{k}}\|_{L^{\infty}(\R^2)}}$, so that $\|\hat{\mathcal{R}}_{a_{k}}\|_{L^{\infty}(\R^2)}=1$. We obtain from (\ref{4.14}), (\ref{4.45}), (\ref{4.50}) and (\ref{4:2}) that
\begin{equation*}
    |\nabla \hat{\mathcal{R}}_{a_{k}}(0)|=\frac{|\nabla \mathrm{R}_{a_{k}}(0)|}{\|\mathcal{R}_{a_{k}}\|_{L^{\infty}(\R^2)}}
    =\frac{1}{\big|Q(0)+\mathrm{R}_{a_{k}}(0)\big|}\frac{|\mathrm{I}_{a_{k}}(0)||\nabla \mathrm{I}_{a_{k}}(0)|}{\|\mathcal{R}_{a_{k}}\|_{L^{\infty}(\R^2)}}
    \leq\frac{C\Omega_{a_{k}}^{2}\bar{\eps}_{a_{k}}^{4}}{\|\mathcal{R}_{a_{k}}\|_{L^{\infty}(\R^2)}}\ \ \hbox{as}\ \ k\to\infty,
\end{equation*}
which then implies that
\begin{equation}\label{4.56}
    \big|\nabla\hat{\mathcal{R}}_{a_{k}}(0)\big|\rightarrow 0\ \ \hbox{as}\ \ k\to\infty.
\end{equation}
It follows from (\ref{4:3}) that
\begin{equation}\label{4.57}
    \tilde{\mathrm{L}}_{a_{k}}\hat{\mathcal{R}}_{a_{k}}(x)
    =\frac{N_{a_{k}}(x)}{\|\mathcal{R}_{a_{k}}\|_{L^{\infty}(\R^2)}}\ \ \hbox{in}\ \ \R^2.
\end{equation}
Note from (\ref{4:7}) that
\begin{equation}\label{4.58}
\begin{split}
    \frac{N_{a_{k}}(x)}{\|\mathcal{R}_{a_{k}}\|_{L^{\infty}(\R^2)}}:&
    =\frac{\Omega_{a_{k}}^{2}\bar{\eps}_{a_{k}}^{4}}{\|\mathcal{R}_{a_{k}}\|_{L^{\infty}(\R^2)}}
    \frac{N_{a_{k}}(x)}{\Omega_{a_{k}}^{2}\bar{\eps}_{a_{k}}^{4}}\\
    &\leq \frac{\Omega_{a_{k}}^{2}\bar{\eps}_{a_{k}}^{4}}{\|\mathcal{R}_{a_{k}}\|_{L^{\infty}(\R^2)}}Ce^{-\frac{1}{10}|x|}\ \ \hbox{uniformly in}\ \ \R^2\ \ \hbox{as}\ \ k\to\infty,
\end{split}
\end{equation}
where the constant $C>0$ is independent of $k>0$. Suppose $y_{a_{k}}$ is a global maximum point of $|\hat{\mathcal{R}}_{a_{k}}(x)|$, so that $|\hat{\mathcal{R}}_{a_{k}}(y_{a_{k}})|=\max\limits_{x\in\R^2}\frac{|\mathcal{R}_{a_{k}}(x)|}{\|\mathcal{R}_{a_{k}}\|_{L^{\infty}(\R^2)}}=1$. By the maximum principle, we deduce from (\ref{4.57}) and (\ref{4.58}) that $|y_{a_{k}}|\leq C$ uniformly in $k>0$.

On the other hand, the elliptic regularity theory shows that there exist a subsequence, still denoted by $\{\hat{\mathcal{R}}_{a_{k}}(x)\}$, of $\{\hat{\mathcal{R}}_{a_{k}}(x)\}$ and a function $\hat{\mathcal{R}}_{0}$ such that $\hat{\mathcal{R}}_{a_{k}}\rightharpoonup \hat{\mathcal{R}}_{0}$ weakly in $H^{1}(\R^2)$ and strongly in $L^{q}_{loc}(\R^2)$ for all $q\in[2,\infty)$ as $k\to\infty$. From (\ref{4.56})--(\ref{4.58}), we get that $\hat{\mathcal{R}}_{0}$ satisfies
\begin{equation*}
    \nabla\hat{\mathcal{R}}_{0}(0)=0,\ \ \tilde{\mathcal{L}}\hat{\mathcal{R}}_{0}(x)=0\ \ \hbox{in}\ \ \R^2,
\end{equation*}
which implies that $\hat{\mathcal{R}}_{0}(x)=\sum\limits_{i=1}^{2}c_{i}\frac{\partial Q}{\partial x_{i}}$ in view of (\ref{4.11}). Since $\nabla\hat{\mathcal{R}}_{0}(0)=0$, we obtain that
\begin{equation}\label{4:8}
    \Big(\frac{\partial^{2}Q(0)}{\partial x_{i}\partial x_{j}}\Big)\left(
     \begin{array}{c}
    c_1\\
    c_2\\
    \end{array}
    \right)=0.
\end{equation}
Because $det\big(\frac{\partial^{2}Q(0)}{\partial x_{i}\partial x_{j}}\big)\neq 0$, we deduce from (\ref{4:8}) that $c_{1}=c_{2}=0$, and hence $\hat{\mathcal{R}}_{0}(x)=0$ in $\R^2$. This contradicts to the fact that up to a subsequence if necessary, $1\equiv\hat{\mathcal{R}}_{a_{k}}(y_{a_{k}})\rightarrow\hat{\mathcal{R}}_{0}(\bar{y}_{0})$ as $k\to\infty$ for some $\bar{y}_{0}\in\R^2$. Therefore, the claim (\ref{4:6}) holds true.

We next claim that
\begin{equation}\label{4:16}
    |\mathcal{R}_{a_{k}}(x)|,\ \ |\nabla\mathcal{R}_{a_{k}}(x)|\leq C\Omega_{a_{k}}^{2}\bar{\eps}_{a_{k}}^{4}e^{-\frac{1}{11}|x|}\ \ \hbox{uniformly in}\ \ \R^2\ \ \hbox{as}\ \ k\to\infty.
\end{equation}
Essentially, we deduce from (\ref{4:3}) and (\ref{4:7}) that
\begin{equation*}
    \tilde{\mathrm{L}}_{a_{k}}\frac{\mathcal{R}_{a_{k}}}{\Omega_{a_{k}}^{2}\bar{\eps}_{a_{k}}^{4}}
    =\frac{N_{a_{k}}(x)}{\Omega_{a_{k}}^{2}\bar{\eps}_{a_{k}}^{4}}
    \leq C_{1}e^{-\frac{1}{10}|x|}\ \ \ \hbox{in}\ \ \R^2\ \ \hbox{as}\ \ k\to\infty.
\end{equation*}
On the other hand, for sufficiently large $R>1$, there exists a constant $C_{2}=C_{2}(R)>0$ such that
\begin{equation}\label{4:10}
    \frac{\mathcal{R}_{a_{k}}}{\Omega_{a_{k}}^{2}\bar{\eps}_{a_{k}}^{4}}\leq C_{2}e^{-\frac{1}{10}|x|}\ \
    \ \hbox{at}\ \ |x|=R\ \ \hbox{as}\ \ k\to\infty,
\end{equation}
and
\begin{equation*}
    C_{2}\tilde{\mathrm{L}}_{a_{k}}e^{-\frac{1}{10}|x|}\geq C_{1}e^{-\frac{1}{10}|x|}\ \ \hbox{in}\ \ \R^2\backslash B_{R}(0)\ \ \hbox{as}\ \ k\to\infty.
\end{equation*}
We thus have
\begin{equation}\label{4:12}
    \tilde{\mathrm{L}}_{a_{k}}\Big(\frac{\mathcal{R}_{a_{k}}}{\Omega_{a_{k}}^{2}\bar{\eps}_{a_{k}}^{4}}
    -C_{2}e^{-\frac{1}{10}|x|}\Big)\leq 0\ \ \ \hbox{in}\ \ \R^2\backslash B_{R}(0)\ \ \hbox{as}\ \ k\to\infty.
\end{equation}
By the comparison principle, we obtain from (\ref{4:10}) and (\ref{4:12}) that
\begin{equation*}
    \frac{\mathcal{R}_{a_{k}}}{\Omega_{a_{k}}^{2}\bar{\eps}_{a_{k}}^{4}}\leq C_{2}e^{-\frac{1}{10}|x|}\ \ \ \hbox{in}\ \ \R^2\backslash B_{R}(0)\ \ \hbox{as}\ \ k\to\infty.
\end{equation*}
Similarly, we deduce that there exists a constant $C_{3}=C_{3}(R)>0$ such that
\begin{equation*}
    \frac{\mathcal{R}_{a_{k}}}{\Omega_{a_{k}}^{2}\bar{\eps}_{a_{k}}^{4}}\geq -C_{3}e^{-\frac{1}{10}|x|}\ \ \ \hbox{in}\ \ \R^2\backslash B_{R}(0)\ \ \hbox{as}\ \ k\to\infty,
\end{equation*}
and therefore
\begin{equation*}
    |\mathcal{R}_{a_{k}}|\leq C\Omega_{a_{k}}^{2}\bar{\eps}_{a_{k}}^{4}e^{-\frac{1}{10}|x|}\ \ \ \hbox{in}\ \ \R^2\backslash B_{R}(0)\ \ \hbox{as}\ \ k\to\infty.
\end{equation*}
This yields that $\mathcal{R}_{a_{k}}$ satisfies
\begin{equation}\label{4:14}
    |\mathcal{R}_{a_{k}}|\leq C\Omega_{a_{k}}^{2}\bar{\eps}_{a_{k}}^{4}e^{-\frac{1}{10}|x|}\ \ \ \hbox{in}\ \ \R^2\ \ \hbox{as}\ \ k\to\infty,
\end{equation}
due to the estimate (\ref{4:6}). By the gradient estimates of (3.15) in \cite{GT}, we further derive from \eqref{4:1}, (\ref{4:3}) and (\ref{4:14}) that
\begin{equation}\label{4:15}
    |\nabla\mathcal{R}_{a_{k}}|\leq C\Omega_{a_{k}}^{2}\bar{\eps}_{a_{k}}^{4}e^{-\frac{1}{11}|x|}\ \ \hbox{in}\ \ \ \R^2\ \ \hbox{as}\ \ k\to\infty.
\end{equation}
We thus obtain from (\ref{4:14}) and (\ref{4:15}) that (\ref{4:16}) holds true.

\emph{Step 2.} Note from (\ref{4.50}) that
\begin{equation*}
    x^{\bot}\cdot\nabla \mathrm{R}_{a_{k}}= x^{\bot}\cdot\nabla \big[\mathcal{R}_{a_{k}}+\Omega_{a_{k}}^{2}\bar{\eps}_{a_{k}}^{4}\Psi_{1}\big]\ \ \hbox{in}\ \ \R^2.
\end{equation*}
Following \eqref{4:1} and (\ref{4:16}), we then obtain that
\begin{equation}\label{4:18}
    |x^{\bot}\cdot\nabla \mathrm{R}_{a_{k}}|\leq C\Omega_{a_{k}}^{2}\bar{\eps}_{a_{k}}^{4}e^{-\frac{1}{12}|x|}\ \ \hbox{uniformly in}\ \ \R^2\ \ \hbox{as}\ \ k\to\infty.
\end{equation}
Similar to (\ref{4:14}) and (\ref{4:15}), we deduce from (\ref{6.9}) and (\ref{4:18}) that
\begin{equation}\label{4:19}
    |\mathrm{I}_{a_{k}}(x)|,\ |\nabla \mathrm{I}_{a_{k}}(x)|\leq C\Omega_{a_{k}}^3\bar{\eps}_{a_{k}}^{6}e^{-\frac{1}{14}|x|}\ \ \hbox{uniformly in}\ \ \R^2\ \ \hbox{as}\ \ k\to\infty,
\end{equation}
where the property (\ref{4.9}) is used in view of the fact that $\int_{\R^2}Q\mathrm{I}_{a_{k}}dx=0$. This yields that (\ref{4.47}) holds true.

\emph{Step 3.} Applying (\ref{4:19}), the same argument of (\ref{4:7}) yields that the nonhomogeneous term $N_{a_{k}}(x)$ of (\ref{4:3}) satisfies
\begin{equation*}
    \frac{|N_{a_{k}}(x)|}{\Omega_{a_{k}}^{2}\bar{\eps}_{a_{k}}^{4}}\leq C\delta_{a_{k}}e^{-\frac{1}{10}|x|}\ \ \hbox{uniformly in}\ \ \R^2\ \ \hbox{as}\ \ k\to\infty,
\end{equation*}
where $\delta_{a_{k}}>0$ satisfies $\delta_{a_{k}}=o(1)$ as $k\to\infty$. The same argument of Step 1 then gives that
\begin{equation}\label{4:21}
    |\mathcal{R}_{a_{k}}|, \ |\nabla \mathcal{R}_{a_{k}}|\leq C\Omega_{a_{k}}^{2}\bar{\eps}_{a_{k}}^{4}\delta_{a_{k}}e^{-\frac{1}{12}|x|}\ \ \hbox{in}\ \ \R^2\ \ \hbox{as}\ \ k\to\infty,
\end{equation}
where $\delta_{a_{k}}>0$ also satisfies $\delta_{a_{k}}=o(1)$ as $k\to\infty$. This proves (\ref{4.48}) in view of (\ref{4.50}),
and Lemma \ref{lem4.4} is therefore proved.
\qed
\vskip 0.05truein

\begin{lem}\label{lem4.5}
Assume $\Omega_{a_{k}}>0$ satisfies \eqref{0:1}, i.e., $\Omega_{a_{k}}=C_{0}(a^*-a_{k})^{-\beta}$ holds for some $\beta\in[0,\frac{1}{2})$ and $C_{0}>0$. Then the imaginary part $\mathrm{I}_{a_{k}}(x)$ of (\ref{6.2}) satisfies
\begin{equation}\label{4.68}
    \mathrm{I}_{a_{k}}(x):=\Omega_{a_{k}}^{3}\bar{\eps}_{a_{k}}^{6}\Phi_{I}(x)
    +o(\Omega_{a_{k}}^{3}\bar{\eps}_{a_{k}}^{6})\ \ \hbox{in}\ \ \R^2\ \ \hbox{as}\ \ k\to\infty,
\end{equation}
where $\Phi_{I}(x)\in C^{2}(\R^2)\cap L^{\infty}(\R^2)$ solves uniquely the following problem
\begin{equation}\label{4.69}
    \int_{\R^2}\Phi_{I}Qdx=0,\ \ \mathcal{L}\Phi_{I}(x)=-2(x^{\bot}\cdot\nabla\Psi_{1})\ \ \hbox{in}\ \ \R^2.
\end{equation}
Here the operator $\mathcal{L}$ is defined by (\ref{4.8}), and the function $\Psi_{1}(x)$ is defined by (\ref{4.49}).
\end{lem}

\noindent\textbf{Proof.} By Lemma \ref{lem4.4} (2), we obtain from (\ref{6.4}) and (\ref{6.9}) that $\int_{\R^2}Q\mathrm{I}_{a_{k}}dx\equiv0$, and
\begin{equation}\label{4.70}
    \mathrm{L}_{a_{k}}\mathrm{I}_{a_{k}}=-2\Omega_{a_{k}}\bar{\eps}_{a_{k}}^{2}(x^{\bot}\cdot\nabla \mathrm{R}_{a_{k}})
    =-2\Omega_{a_{k}}^{3}\bar{\eps}_{a_{k}}^{6}(x^{\bot}\cdot\nabla \Psi_{1})+o(\Omega_{a_{k}}^{3}\bar{\eps}_{a_{k}}^{6})\ \ \hbox{as}\ \ k\to\infty,
\end{equation}
where $\Psi_{1}(x)$ is defined by (\ref{4.49}). Set
\begin{equation*}
    \mathrm{I}_{1,a_{k}}(x)=\mathrm{I}_{a_{k}}(x)-\Omega_{a_{k}}^{3}\bar{\eps}_{a_{k}}^{6}\Phi_{I}(x),
\end{equation*}
where $\Phi_{I}(x)\in C^{2}(\R^2)\cap L^{\infty}(\R^2)$ is defined by (\ref{4.69}).
Similar to Step 3 in the proof of Lemma \ref{lem4.4}, we then deduce from (\ref{4.47}) and (\ref{4.70}) that $\mathrm{I}_{1,a_{k}}(x)=o(\Omega_{a_{k}}^{3}\bar{\eps}_{a_{k}}^{6})$ uniformly in $\R^2$ as $k\to\infty$, which implies that (\ref{4.68}) holds true. Also, the restriction $\int_{\R^2}\Phi_{I}Qdx=0$ and the property (\ref{4.9}) further give the uniqueness of $\Phi_{I}(x)\in C^{2}(\R^2)\cap L^{\infty}(\R^2)$ defined by (\ref{4.69}), and we are hence done.
\qed

\section{Nonexistence of Vortices}
In this section, we shall first establish the refined estimates of the unique global maximum point $x_{a_{k}}\in\R^2$ obtained in \eqref{4:22} for $|u_{a_{k}}|$ as $k\to\infty$, by which we then prove Theorem \ref{th5} on the refined expansions of $u_{a_{k}}$ as $k\to\infty$. We finally prove Theorem \ref{th4} by applying Theorem \ref{th5}.

\begin{lem}\label{lem6.4}
Let $u_{a_{k}}$ be a complex-valued minimizer of $I(a_{k},\Omega_{a_{k}})$, where $\Omega_{a_{k}}>0$ satisfies \eqref{0:1} for some $\beta\in[0,\frac{1}{2})$ and $C_{0}>0$. Let $x_{a_{k}}\in\R^2$ be the unique maximum point of $|u_{a_{k}}(x)|$ satisfying $\lim\limits_{k\to\infty} x_{a_{k}}=x_{0}$, where $|x_0|=1$. Then we have
\begin{equation}\label{5.02}
    |x_{a_{k}}|^2-1=\tilde{C}\bar{\eps}_{a_{k}}^2+o(\Omega_{a_{k}}\bar{\eps}_{a_{k}}^4)\ \ \hbox{as}\ \ k\to\infty,
\end{equation}
where $\bar{\eps}_{a_{k}}>0$ is as in \eqref{6.1}, and the constant $\tilde{C}=-\frac{8\lambda^4}{5a^*}<0$.
\end{lem}

\noindent\textbf{Proof.}
We first consider the case where $x_{0}\neq (0,1)$. It then follows from \eqref{6.10} and \eqref{6.11} that
\begin{equation}\label{6.029}
\begin{aligned}
    \int_{\R^2}\frac{\partial Q}{\partial x_{1}}\tilde{\mathrm{L}}_{a_{k}}\mathrm{R}_{a_{k}}dx
    &=\int_{\R^2}\frac{\partial Q}{\partial x_{1}}\mathrm{F}_{a_{k}}(x)dx\\
    &=\int_{\R^2}\frac{\partial Q}{\partial x_{1}}\Big\{-\Big[\Omega_{a_{k}}^{2}\bar{\varepsilon}_{a_{k}}^{4}|x|^2
    +\frac{\Omega_{a_{k}}^{2}\bar{\varepsilon}_{a_{k}}^{2}}{8}(|\bar{\varepsilon}_{a_{k}}x+x_{a_{k}}|^{2}-1)^2\Big]Q\\
   &\qquad+2\Omega_{a_{k}}\bar{\varepsilon}_{a_{k}}^{2}(x^{\bot}\cdot\nabla \mathrm{I}_{a_{k}})+\mathrm{I}_{a_{k}}^{2}Q\Big\}dx\\
    &=-\frac{1}{4}\Omega_{a_{k}}^2\bar{\eps}_{a_{k}}^4\int_{\R^2}\frac{\partial Q^2}{\partial x_{1}}\Big(\bar{\eps}_{a_{k}}|x|^2+\frac{|x_{a_{k}}|^2-1}{\bar{\eps}_{a_{k}}}\Big)(x\cdot x_{a_{k}})dx\\
   &\quad+2\Omega_{a_{k}}\bar{\eps}_{a_{k}}^2\int_{\R^2}\frac{\partial Q}{\partial x_{1}}(x^{\bot}\cdot\nabla\mathrm{I}_{a_{k}})dx
    +\int_{\R^2}\frac{\partial Q}{\partial x_{1}}\mathrm{I}_{a_{k}}^2Qdx,
\end{aligned}
\end{equation}
where the operator $\tilde{\mathrm{L}}_{a_{k}}$ is defined by \eqref{6.10}. We note that
\begin{equation}\label{6.030}
\begin{aligned}
&2\Omega_{a_{k}}\bar{\eps}_{a_{k}}^2\int_{\R^2}\frac{\partial Q}{\partial x_{1}}(x^{\bot}\cdot\nabla\mathrm{I}_{a_{k}})dx
    +\int_{\R^2}\frac{\partial Q}{\partial x_{1}}\mathrm{I}_{a_{k}}^2Qdx\\
    =&\big(2\Omega_{a_{k}}^4\bar{\eps}_{a_{k}}^8
    +o(\Omega_{a_{k}}^4\bar{\eps}_{a_{k}}^8)\big)\int_{\R^2}\frac{\partial Q}{\partial x_{1}}(x^{\bot}\cdot\nabla\Phi_{I})dx\ \ \hbox{as}\ \ k\to\infty
\end{aligned}
\end{equation}
in view of \eqref{4.68}, where $\Phi_{I}$ is defined by \eqref{4.69}.
We thus obtain from \eqref{6.029} and \eqref{6.030} that
\begin{equation}\label{6.031}
\begin{aligned}
    &\quad-\frac{1}{4}\Omega_{a_{k}}^2\bar{\eps}_{a_{k}}^3(|x_{a_{k}}|^2-1)\int_{\R^2}\frac{\partial Q^2}{\partial x_{1}}(x\cdot x_{a_{k}})dx\\
    &=\frac{1}{4}\Omega_{a_{k}}^2\bar{\eps}_{a_{k}}^5\int_{\R^2}\frac{\partial Q^2}{\partial x_{1}}|x|^2(x\cdot x_{a_{k}})dx
    +\int_{\R^2}\frac{\partial Q}{\partial x_{1}}\tilde{\mathrm{L}}_{a_{k}}\mathrm{R}_{a_{k}}dx\\
    &\quad-2\Omega_{a_{k}}^4\bar{\eps}_{a_{k}}^8\int_{\R^2}\frac{\partial Q}{\partial x_{1}}\big(x^{\bot}\cdot\nabla\Phi_{I}\big)dx+o(\Omega_{a_{k}}^4\bar{\eps}_{a_{k}}^8)
    \ \ \hbox{as}\ \ k\to\infty.
\end{aligned}
\end{equation}

We next simplify the right hand side of \eqref{6.031}. By the definition of $\tilde{\mathcal{L}}$ in (\ref{4.10}), we have $\int_{\R^2}\frac{\partial Q}{\partial x_{1}}\tilde{\mathcal{L}}\mathrm{R}_{a_{k}}dx=0$. We then deduce from \eqref{6.10} and \eqref{4.48} that
\begin{equation}\label{6.032}
\begin{split}
   &\quad \int_{\R^2}\frac{\partial Q}{\partial x_{1}}\tilde{\mathrm{L}}_{a_{k}}\mathrm{R}_{a_{k}}dx\\
    &=\int_{\R^2}\frac{\partial Q}{\partial x_{1}}(\tilde{\mathrm{L}}_{a_{k}}-\tilde{\mathcal{L}})\mathrm{R}_{a_{k}}dx\\
    &=\int_{\R^2}\frac{\partial Q}{\partial x_{1}}\Big\{\Omega_{a_{k}}^2\bar{\eps}_{a_{k}}^4|x|^2
    +\frac{\Omega_{a_{k}}^2\bar{\eps}_{a_{k}}^2}{8}(|\bar{\eps}_{a_{k}}x+x_{a_{k}}|^2-1)^2
    -(3Q+\mathrm{R}_{a_{k}})\mathrm{R}_{a_{k}}-\mathrm{I}_{a_{k}}^2\Big\}\mathrm{R}_{a_{k}}dx\\
    &=\int_{\R^2}\frac{\partial Q}{\partial x_{1}}\Big\{\Omega_{a_{k}}^2\bar{\eps}_{a_{k}}^4
    \Big[|x|^2+\frac{1}{8}\Big(\bar{\eps}_{a_{k}}^2|x|^4+4(x\cdot x_{a_{k}})^2+\Big(\frac{|x_{a_{k}}|^2-1}{\bar{\eps}_{a_{k}}}\Big)^2\\
    &\quad+4\bar{\eps}_{a_{k}}|x|^2(x\cdot x_{a_{k}})+2(|x_{a_{k}}|^2-1)|x|^2+4\frac{|x_{a_{k}}|^2-1}{\bar{\eps}_{a_{k}}}(x\cdot x_{a_{k}})\Big)\Big]\\
    &\quad-(3Q+\mathrm{R}_{a_{k}})\mathrm{R}_{a_{k}}-\mathrm{I}_{a_{k}}^2\Big\}\mathrm{R}_{a_{k}}dx\\
    &=\Omega_{a_{k}}^4\bar{\eps}_{a_{k}}^8\Big\{\int_{\R^2}\frac{\partial Q}{\partial x_{1}}\big[|x|^2+\frac{1}{2}(x\cdot x_{0})^2\big]\Psi_{1}dx-3\int_{\R^2}\frac{\partial Q}{\partial x_{1}}Q\Psi_{1}^2dx\Big\}+o(\Omega_{a_{k}}^4\bar{\eps}_{a_{k}}^8)\\
\end{split}
\end{equation}
as $k\to\infty$,
where $\Psi_{1}$ is defined by \eqref{4.49}.
We conclude from \eqref{6.031} and \eqref{6.032} that
\begin{equation}\label{6.033}
\begin{aligned}
    &-\frac{1}{4}\Omega_{a_{k}}^2\bar{\eps}_{a_{k}}^3(|x_{a_{k}}|^2-1)\int_{\R^2}\frac{\partial Q^2}{\partial x_{1}}(x\cdot x_{a_{k}})dx\\
    &=\frac{1}{4}\Omega_{a_{k}}^2\bar{\eps}_{a_{k}}^5\int_{\R^2}\frac{\partial Q^2}{\partial x_{1}}|x|^2(x\cdot x_{a_{k}})dx
    +\Omega_{a_{k}}^4\bar{\eps}_{a_{k}}^8\Big\{\int_{\R^2}\frac{\partial Q}{\partial x_{1}}\Big[|x|^2+\frac{1}{2}(x\cdot x_{0})^2\Big]\Psi_{1}dx\\
    &\quad-3\int_{\R^2}\frac{\partial Q}{\partial x_{1}}Q\Psi_{1}^2dx-2\int_{\R^2}\frac{\partial Q}{\partial x_{1}}(x^{\bot}\cdot\nabla\Phi_{I})dx\Big\}+o(\Omega_{a_{k}}^4\bar{\eps}_{a_{k}}^8)\ \ \hbox{as}\ \ k\to\infty,
\end{aligned}
\end{equation}
where $\Omega_{a_{k}}=C_{0}(a^*-a_{k})^{-\beta}$ holds for some $\beta\in [0,\frac{1}{2})$ and $C_{0}>0$, and $\bar{\eps}_{a_{k}}=\frac{(a^*-a_{k})^{\frac{1+2\beta}{4}}}{\sqrt{C_{0}}\lambda}\big[1+o(1)\big]$ as $k\to\infty$.

Since $x_{0}\neq (0,1)$, we have
\begin{equation*}
    \int_{\R^2}\frac{\partial Q^2}{\partial x_{1}}(x\cdot x_{0})dx\neq 0,
\end{equation*}
and we hence derive from \eqref{6.033} that
\begin{equation}\label{6:1}
    \big||x_{a_{k}}|^2-1\big|\leq C\bar{\eps}_{a_{k}}^2\ \ \mbox{as}\ \ k\to\infty.
\end{equation}
Inserting (\ref{6:1}) into (\ref{6.033}) yields  that
\begin{equation}\label{5:1}
    |x_{a_{k}}|^2-1=\bar{\eps}_{a_{k}}^2\tilde{C}+o(\Omega_{a_{k}}\bar{\eps}_{a_{k}}^4)\ \ \hbox{as}\ \ k\to\infty.
\end{equation}
The constant $\tilde{C}$ in \eqref{5:1} satisfies
\begin{equation*}
\begin{aligned}
  \tilde{C}=-\frac{\int_{\R^2}\frac{\partial Q^2}{\partial x_{1}}|x|^2(x\cdot x_{0})dx}{\int_{\R^2}\frac{\partial Q^2}{\partial x_{1}}(x\cdot x_{0})dx}&=-\frac{x_{01}\int_{\R^2}(2x_{1}^2+|x|^2)Q^2(x)dx}{x_{01}\int_{\R^2}Q^2(x)dx}\\
  &=-\frac{2\int_{\R^2}|x|^2Q^2dx}{a^*}
  =-\frac{8\lambda^4}{5a^*},
\end{aligned}
\end{equation*}
where $x_{0}=(x_{01},x_{02})\in\R^2$, $\lambda>0$ is as in \eqref{4.1}, and $Q(x)=Q(|x|)$ is also used.
The lemma is therefore proved for the   case where $x_{0}\not =(0,1)$.

We next consider the case where $x_{0}=(0,1)$. Similar to \eqref{6.029}--\eqref{6.033}, we then obtain that
\begin{equation*}
\begin{aligned}
    &-\frac{1}{4}\Omega_{a_{k}}^2\bar{\eps}_{a_{k}}^3(|x_{a_{k}}|^2-1)\int_{\R^2}\frac{\partial Q^2}{\partial x_{2}}(x\cdot x_{a_{k}})dx\\
    &=\frac{1}{4}\Omega_{a_{k}}^2\bar{\eps}_{a_{k}}^5\int_{\R^2}\frac{\partial Q^2}{\partial x_{2}}|x|^2(x\cdot x_{a_{k}})dx
    +\Omega_{a_{k}}^4\bar{\eps}_{a_{k}}^8\Big\{\int_{\R^2}\frac{\partial Q}{\partial x_{2}}\Big[|x|^2+\frac{1}{2}(x\cdot x_{0})^2\Big]\Psi_{1}dx\\
    &\quad-3\int_{\R^2}\frac{\partial Q}{\partial x_{2}}Q\Psi_{1}^2dx-2\int_{\R^2}\frac{\partial Q}{\partial x_{2}}(x^{\bot}\cdot\nabla\Phi_{I})dx\Big\}+o(\Omega_{a_{k}}^4\bar{\eps}_{a_{k}}^8)\ \ \hbox{as}\ \ k\to\infty,
\end{aligned}
\end{equation*}
where as before $\Omega_{a_{k}}=C_{0}(a^*-a_{k})^{-\beta}$ holds for some $\beta\in [0,\frac{1}{2})$ and $C_{0}>0$, and $\bar{\eps}_{a_{k}}=\frac{(a^*-a_{k})^{\frac{1+2\beta}{4}}}{\sqrt{C_{0}}\lambda}[1+o(1)]$ as $k\to\infty$.
Since $x_{0}=(0,1)$, we have
\begin{equation*}
    \int_{\R^2}\frac{\partial Q^2}{\partial x_{2}}(x\cdot x_{0})dx\neq 0.
\end{equation*}
The same arguments of \eqref{6:1}--\eqref{5:1} then yield that (\ref{5.02}) holds for
\begin{equation*}
  \tilde{C}=-\frac{\int_{\R^2}\frac{\partial Q^2}{\partial x_{2}}|x|^2(x\cdot x_{0})dx}{\int_{\R^2}\frac{\partial Q^2}{\partial x_{2}}(x\cdot x_{0})dx}
  =-\frac{\int_{\R^2}(2x_{2}^2+|x|^2)Q^2(x)dx}{\int_{\R^2}Q^2(x)dx}
  =-\frac{8\lambda^4}{5a^*},
\end{equation*}
where $\lambda>0$ is defined by \eqref{4.1}.
This completes the proof of Lemma \ref{lem6.4}.
\qed
\vskip 0.05truein

Following Lemma \ref{lem6.4}, we are now ready to derive the following refined expansion of $u_{a_{k}}(x)$ as $k\to\infty$.

\begin{thm}\label{th5}
Assume $u_{a_{k}}$ is a complex-valued minimizer of $I(a_{k},\Omega_{a_{k}})$, where $\Omega_{a_{k}}>0$ satisfies \eqref{0:1} for some $\beta\in[0,\frac{1}{2})$ and $C_{0}>0$. Then we have
\begin{equation}\label{K:E1}
\begin{split}
    \nu_{a_{k}}(x):
    &=\bar{\eps}_{a_{k}}\sqrt{a_{k}}u_{a_{k}}(\bar{\eps}_{a_{k}}x+x_{a_{k}})e^{-i(\Omega_{a_{k}}\bar{\eps}_{a_{k}}x\cdot x_{a_{k}}^{\bot}-\rho_{a_{k}})}\\
    &=Q(x)+\Omega_{a_{k}}^2\bar{\eps}_{a_{k}}^4\Psi_{1}(x)
    +\Omega_{a_{k}}^2\bar{\eps}_{a_{k}}^5\Psi_{2}(x)\\
    &\quad+i\Omega_{a_{k}}^3\bar{\eps}_{a_{k}}^6\Phi_{I}(x)+o(\Omega_{a_{k}}^3\bar{\eps}_{a_{k}}^6)
    \ \ \hbox{as}\ \ k\to\infty,
\end{split}
\end{equation}
where $\bar{\eps}_{a_{k}}>0$ is as in \eqref{6.1}, $\rho_{a_{k}}\in [0,2\pi)$ is a suitable constant phase, and $x_{a_{k}}$ is the unique global maximum point of $|u_{a_{k}}|$ and satisfies
\begin{equation*}
 |x_{a_{k}}|^2-1=\tilde{C}\bar{\eps}_{a_{k}}^2+o(\Omega_{a_{k}}\bar{\eps}_{a_{k}}^4)\ \ \hbox{as}\ \ k\to\infty.
\end{equation*}
Here $\Psi_{i}$ is defined by \eqref{0.1} for $i=1,2$, and $\Phi_{I}$ is defined by \eqref{0.3}.
\end{thm}

\noindent\textbf{Proof.}
Under the assumptions of Theorem \ref{th5}, similar to Lemma \ref{lem4.4}, one can deduce from Lemmas \ref{lem4.4} and \ref{lem6.4} that the real part of \eqref{6.2} satisfies
\begin{equation*}
    \tilde{\mathrm{R}}_{a_{k}}(x):=Q(x)+\Omega_{a_{k}}^2\bar{\eps}_{a_{k}}^4\Psi_{1}(x)
    +\Omega_{a_{k}}^2\bar{\eps}_{a_{k}}^5\Psi_{2}(x)+O(\Omega_{a_{k}}^2\bar{\eps}_{a_{k}}^6)\ \ \hbox{as}\ \ k\to\infty,
\end{equation*}
where $\Psi_{1}$ is defined by \eqref{4.49}, and $\Psi_{2}(x)\in C^2(\R^2)\bigcap L^{\infty}(\R^2)$ solves uniquely the following problem
\begin{equation*}
    \nabla\Psi_{2}(0)=0,\ \tilde{\mathcal{L}}\Psi_{2}(x)=-\frac{1}{2}\big(|x|^2+\tilde{C}\big)(x\cdot x_{0})Q(x).
\end{equation*}
Here the constant $\tilde{C}=-\frac{8\lambda^4}{5a^*}<0$. The proof of Theorem \ref{th5} is therefore complete in view of Lemma \ref{lem4.5}.
\qed
\vskip 0.05truein

\subsection{Proof of Theorem \ref{th4}}

Applying Theorem \ref{th5}, we finally address the proof of Theorem \ref{th4} on the non-existence of vortices in a very large region.
\vskip 0.05truein
\noindent\textbf{Proof of Theorem \ref{th4}.}
It follows from Theorem \ref{th5} that for any fixed and sufficiently large $R>1$, there exists a constant $C_{R}:=C(R)>0$ such that
\begin{equation}\label{6.38}
    |\nu_{a_{k}}(x)|\geq Q(x)-C_{R}\Omega_{a_{k}}^2\bar{\eps}_{a_{k}}^4>0\ \ \hbox{in}\ \ \{x\in\R^2:|x|\leq R\}\ \ \hbox{as}\ \ k\to\infty,
\end{equation}
where $\nu_{a_{k}}(x)$ is as in (\ref{K:E1}), $\Omega_{a_{k}}=C_{0}(a^*-a_{k})^{-\beta}$ holds for some $\beta\in [0,\frac{1}{6})$ and $C_{0}>0$, and $\bar{\eps}_{a_{k}}=\frac{(a^*-a_{k})^{\frac{1+2\beta}{4}}}{\sqrt{C_{0}}\lambda}\big[1+o(1)\big]$ as $k\to\infty$.

Denoting $\tilde{w}_{a_{k}}=\nu_{a_{k}}-Q(x)$, we deduce from (\ref{6.7}) that $\tilde{w}_{a_{k}}$ satisfies
\begin{equation*}
    (-\Delta+\hat{V}_{a_{k}})(\tilde{w}_{a_{k}}+Q)+2i\Omega_{a_{k}}\bar{\eps}_{a_{k}}^{2} (x^{\bot}\cdot\nabla\tilde{w}_{a_{k}})=0
    \ \ \hbox{in}\ \ \R^2,
\end{equation*}
i.e.,
\begin{equation}\label{1}
    (-\Delta+\hat{V}_{a_{k}})\tilde{w}_{a_{k}}+2i\Omega_{a_{k}}\bar{\eps}_{a_{k}}^{2} (x^{\bot}\cdot\nabla\tilde{w}_{a_{k}})+(-\Delta+\hat{V}_{a_{k}})Q=0\ \ \hbox{in}\ \ \R^2,
\end{equation}
where
\begin{equation}\label{2}
    \hat{V}_{a_{k}}(x):=\Omega_{a_{k}}^{2}\bar{\eps}_{a_{k}}^{4}|x|^2
    +\frac{\Omega_{a_{k}}^{2}\bar{\eps}_{a_{k}}^{2}}{8}\big(|\bar{\eps}_{a_{k}}x+x_{a_{k}}|^{2}-1\big)^{2}
    +1-|\nu_{a_{k}}|^{2}\ \ \hbox{in}\ \ \R^2.
\end{equation}
Applying (\ref{1}) and (\ref{2}), we have
\begin{equation}\label{3}
\begin{split}
    &-\frac{1}{2}\Delta |\tilde{w}_{a_{k}}|^{2}+\Big[\Omega_{a_{k}}^{2}\bar{\eps}_{a_{k}}^{4}|x|^2
    +\frac{\Omega_{a_{k}}^{2}\bar{\eps}_{a_{k}}^{2}}{8}(|\bar{\eps}_{a_{k}}x+x_{a_{k}}|^{2}-1)^{2}
    +1-|\nu_{a_{k}}|^{2}\Big]|\tilde{w}_{a_{k}}|^{2}    \\
    &+|\nabla\tilde{w}_{a_{k}}|^{2}-2\Omega_{a_{k}} \bar{\eps}_{a_{k}}^{2} x^{\bot}\cdot(i\tilde{w}_{a_{k}},\nabla\tilde{w}_{a_{k}})
    +(-\Delta+\hat{V}_{a_{k}})(Q,\tilde{w}_{a_{k}})=0\ \ \hbox{in}\ \ \R^2.
\end{split}
\end{equation}
Employing the diamagnetic inequality, we obtain that
\begin{equation}\label{4}
    |\nabla\tilde{w}_{a_{k}}|^{2}+\Omega_{a_{k}}^{2}\bar{\eps}_{a_{k}}^{4}|x|^2|\tilde{w}_{a_{k}}|^{2}
    -2\Omega_{a_{k}}\bar{\eps}_{a_{k}}^{2} x^{\bot}\cdot(i\tilde{w}_{a_{k}},\nabla\tilde{w}_{a_{k}})
    \geq \big|\nabla|\tilde{w}_{a_{k}}|\big|^{2}\ \ \ \hbox{in}
    \ \ \R^2.
\end{equation}
We then derive from (\ref{3}) and (\ref{4}) that
\begin{equation}\label{5}
    -\Delta|\tilde{w}_{a_{k}}|+|\tilde{w}_{a_{k}}|
    \leq \big|(-\Delta+\hat{V}_{a_{k}})Q\big|+|\nu_{a_{k}}|^{2}|\tilde{w}_{a_{k}}|\ \ \hbox{in}\ \ \R^2,
\end{equation}
due to the fact that
\[
- \frac{1}{2}\Delta  |\tilde{w}_{a_{k}}|^{2}+\big|\nabla |\tilde{w}_{a_{k}}|\big|^2=-|\tilde{w}_{a_{k}}|\,\Delta |\tilde{w}_{a_{k}}|.
\]

The same argument of \cite[Proposition 3.3]{GLY} yields that there exists a constant $C>0$ such that
\begin{equation}\label{5.05}
    |\nu_{a_{k}}|^2\leq Ce^{-\frac{4}{3}|x|}\ \ \hbox{in}\ \ \R^2\ \ \hbox{as}\ \ k\to\infty,
\end{equation}
together with Theorem \ref{th5}, which further implies  that
\begin{equation}\label{6}
    |\nu_{a_{k}}|^{2}|\tilde{w}_{a_{k}}|\leq C\Omega_{a_{k}}^{2}\bar{\eps}_{a_{k}}^{4}e^{-\frac{4}{3}|x|}\ \ \hbox{in}\ \ \R^2\ \ \hbox{as}\ \ k\to\infty.
\end{equation}
Following (\ref{expo}), \eqref{5.02} and (\ref{2}), we derive that
\begin{equation}\label{7}
\begin{split}
   &\big|(-\Delta+\hat{V}_{a_{k}})Q\big|\\
   =&\Big|\big[\Omega_{a_{k}}^{2}\bar{\eps}_{a_{k}}^{4}|x|^2
    +\frac{\Omega_{a_{k}}^{2}\bar{\eps}_{a_{k}}^{2}}{8}(|\bar{\eps}_{a_{k}}x+x_{a_{k}}|^{2}-1)^{2}
    +Q^{2}-|\nu_{a_{k}}|^{2}\big]Q\Big|\\
    \leq &C \big(\Omega_{a_{k}}^2\bar{\varepsilon}_{a_{k}}^4|x|^\frac{3}{2}
    +\Omega_{a_{k}}^2\bar{\varepsilon}_{a_{k}}^5|x|^\frac{5}{2}
    +\Omega_{a_{k}}^2\bar{\varepsilon}_{a_{k}}^6|x|^{\frac{7}{2}}\big)e^{-|x|}\ \ \hbox{in}\ \ \R^2\backslash B_{R}(0)\ \ \hbox{as}\ \ k\to\infty,
\end{split}
\end{equation}
where the sufficiently large constant $R>1$ is as in \eqref{6.38}.
We then obtain from (\ref{5})--(\ref{7}) that
\begin{equation}\label{A1}
\begin{split}
    -\Delta|\tilde{w}_{a_{k}}|+|\tilde{w}_{a_{k}}|
    \leq &\bar{C}_{0} \big(\Omega_{a_{k}}^2\bar{\varepsilon}_{a_{k}}^4|x|^\frac{3}{2}
    +\Omega_{a_{k}}^2\bar{\varepsilon}_{a_{k}}^5|x|^\frac{5}{2}
    +\Omega_{a_{k}}^2\bar{\varepsilon}_{a_{k}}^6|x|^{\frac{7}{2}}\big)e^{-|x|}\\
        \ \ &\hbox{in}\ \ \R^2\backslash B_{R}(0)\ \ \hbox{as}\ \ k\to\infty,
\end{split}
\end{equation}
where $\bar{C}_{0}>0$ is a constant independent of $k>0$.
Because Theorem \ref{th5}
implies that $|\tilde{w}_{a_{k}}|= O\big(\Omega_{a_{k}}^2\bar{\varepsilon}_{a_{k}}^4\big)$ as $k\to\infty$, we deduce that for the sufficiently large $R>1$ in \eqref{6.38},
\begin{equation*}
  |\tilde{w}_{a_{k}}|\leq  \bar C\big(\Omega_{a_{k}}^2\bar{\varepsilon}_{a_{k}}^4|x|^\frac{5}{2}
    +\Omega_{a_{k}}^2\bar{\varepsilon}_{a_{k}}^5|x|^\frac{7}{2}
    +\Omega_{a_{k}}^2\bar{\varepsilon}_{a_{k}}^6|x|^{\frac{9}{2}}\big)e^{-|x|} \,\ \mbox{at}\, \ |x|=R>1 \,\ \mbox{as}\, \    k\to\infty,
\end{equation*}
where the constant $ \bar C>0$ is independent of $k>0$.

Consider
\begin{equation}\label{6.003}
    \tilde{w}_{1a_{k}}:=C_{1}\Omega_{a_{k}}^2\bar{\varepsilon}_{a_{k}}^4|x|^{\frac{5}{2}}e^{-|x|}\ \ \hbox{in}\ \ \R^2,
\end{equation}
where the constant $C_{1}>0$ is chosen such that as $k\to\infty$,
\begin{equation*}
    -\Delta\tilde{w}_{1a_{k}}+\tilde{w}_{1a_{k}}\geq \bar{C}_{0}\Omega_{a_{k}}^2\bar{\varepsilon}_{a_{k}}^4|x|^{\frac{3}{2}}e^{-|x|}\ \ \hbox{in}\ \ \R^2\backslash B_{R}(0),
\end{equation*}
and
\begin{equation*}
    \tilde{w}_{1a_{k}}\geq\bar{C}\Omega_{a_{k}}^2\bar{\varepsilon}_{a_{k}}^4|x|^\frac{5}{2}e^{-|x|}\ \ \hbox{at}\ \ |x|=R.
\end{equation*}
Similarly, we choose
\begin{equation*}
    \tilde{w}_{2a_{k}}:=C_{2}\Omega_{a_{k}}^2\bar{\varepsilon}_{a_{k}}^5|x|^{\frac{7}{2}}e^{-|x|}
    \ \ \hbox{and}\ \ \tilde{w}_{3a_{k}}:=C_{3}\Omega_{a_{k}}^2\bar{\varepsilon}_{a_{k}}^6|x|^{\frac{9}{2}}e^{-|x|}\ \ \hbox{in}\ \ \R^2,
\end{equation*}
where the constants $C_{2}$ and $C_{3}>0$ are chosen such that as $k\to\infty$,
\begin{equation*}
    -\Delta\tilde{w}_{2a_{k}}+\tilde{w}_{2a_{k}}\geq \bar{C}_{0}\Omega_{a_{k}}^2\bar{\varepsilon}_{a_{k}}^5|x|^{\frac{5}{2}}e^{-|x|},\ -\Delta\tilde{w}_{3a_{k}}+\tilde{w}_{3a_{k}}\geq \bar{C}_{0}\Omega_{a_{k}}^2\bar{\varepsilon}_{a_{k}}^6|x|^{\frac{7}{2}}e^{-|x|}\ \ \hbox{in}\ \ \R^2\backslash B_{R}(0),
\end{equation*}
and
\begin{equation*}
    \tilde{w}_{2a_{k}}\geq\bar{C}\Omega_{a_{k}}^2\bar{\varepsilon}_{a_{k}}^5|x|^\frac{7}{2}e^{-|x|},\ \tilde{w}_{3a_{k}}\geq\bar{C}\Omega_{a_{k}}^2\bar{\varepsilon}_{a_{k}}^6|x|^\frac{9}{2}e^{-|x|}\ \ \hbox{at}\ \ |x|=R.
\end{equation*}
We then derive that as $k\to\infty$,
\begin{equation*}
\begin{aligned}
  &-\Delta(\tilde{w}_{1a_{k}}+\tilde{w}_{2a_{k}}+\tilde{w}_{3a_{k}})
  +(\tilde{w}_{1a_{k}}+\tilde{w}_{2a_{k}}+\tilde{w}_{3a_{k}})\\
  \geq&
  \bar{C}_{0} \big(\Omega_{a_{k}}^2\bar{\varepsilon}_{a_{k}}^4|x|^\frac{3}{2}
    +\Omega_{a_{k}}^2\bar{\varepsilon}_{a_{k}}^5|x|^\frac{5}{2}
    +\Omega_{a_{k}}^2\bar{\varepsilon}_{a_{k}}^6|x|^{\frac{7}{2}}\big)e^{-|x|}
    \ \ \hbox{in}\ \ \R^2\backslash B_{R}(0),
  \end{aligned}
\end{equation*}
and
\begin{equation*}
\begin{aligned}
  \tilde{w}_{1a_{k}}+\tilde{w}_{2a_{k}}+\tilde{w}_{3a_{k}}
  \geq
  \bar{C}(\Omega_{a_{k}}^2\bar{\varepsilon}_{a_{k}}^4|x|^\frac{5}{2}
    +\Omega_{a_{k}}^2\bar{\varepsilon}_{a_{k}}^5|x|^\frac{7}{2}
    +\Omega_{a_{k}}^2\bar{\varepsilon}_{a_{k}}^6|x|^{\frac{9}{2}}\big)e^{-|x|}\ge |\tilde{w}_{a_{k}}|\ \ \hbox{at}\ \ |x|=R,
  \end{aligned}
\end{equation*}
provided that the constant $R>1$ is sufficiently large.
By the comparison principle, we conclude from above that
\begin{equation}\label{6.004}
    |\tilde{w}_{a_{k}}|\leq \tilde{w}_{1a_{k}}+\tilde{w}_{2a_{k}}+\tilde{w}_{3a_{k}}\ \ \hbox{in}\ \ \R^2\backslash B_{R}(0)\ \ \hbox{as}\ \ k\to\infty,
\end{equation}
provided that the constant $R>1$ is sufficiently large.

By the definition of $\tilde{w}_{1a_{k}}$ in \eqref{6.003}, we derive that there exists a constant $C_{0}^{'}>0$ such that as $k\to\infty$,
\begin{equation}\label{6.005}
\begin{aligned}
    \frac{1}{3}Q-\tilde{w}_{1a_{k}}&\geq \frac{1}{3}C_{0}^{'}|x|^{-\frac{1}{2}}e^{-|x|}-C_{1}\Omega_{a_{k}}^2\bar{\varepsilon}_{a_{k}}^4|x|^{\frac{5}{2}}e^{-|x|}\\
    &=|x|^{-\frac{1}{2}}e^{-|x|}\Big(\frac{1}{3}C_{0}^{'}-C_{1}\Omega_{a_{k}}^2\bar{\varepsilon}_{a_{k}}^4|x|^3\Big)>0,
\end{aligned}
\end{equation}
if $R\leq|x|<\big(\frac{C_{0}^{'}}{3C_{1}}\big)^{\frac{1}{3}}(\Omega_{a_{k}}^2\bar{\varepsilon}_{a_{k}}^4)^{-\frac{1}{3}}$. Similar to \eqref{6.005}, we also obtain that
\begin{equation}\label{6.008}
    \frac{1}{3}Q-\tilde{w}_{2a_{k}}>0,\ \ \hbox{if}\ \ R\leq|x|< (\frac{C_{0}^{'}}{3C_{2}})^{\frac{1}{4}}(\Omega_{a_{k}}^2\bar{\eps}_{a_{k}}^5)^{-\frac{1}{4}}\ \ \hbox{as}\ \  k\to\infty;
\end{equation}
and
\begin{equation}\label{6.006}
    \frac{1}{3}Q-\tilde{w}_{3a_{k}}>0,\ \ \hbox{if}\ \ R\leq|x|< (\frac{C_{0}^{'}}{3C_{3}})^{\frac{1}{5}}(\Omega_{a_{k}}^2\bar{\eps}_{a_{k}}^6)^{-\frac{1}{5}}\ \ \hbox{as}\ \  k\to\infty.
\end{equation}
Since $\bar{\eps}_{a_{k}}=\eps_{a_{k}}[1+o(1)]$ as $k\to\infty$, it follows from \eqref{6.004}--\eqref{6.006} that there exists a constant $C_{4}>0$ such that
\begin{equation*}
\begin{aligned}
    |\nu_{a_{k}}|\geq Q-|\nu_{a_{k}}-Q|=Q-|\tilde{w}_{a_{k}}|&\geq Q-\big(\tilde{w}_{1a_{k}}+\tilde{w}_{2a_{k}}+\tilde{w}_{3a_{k}}\big)\\
    &>0,\ \ \hbox{if}\ \    R\leq|x|\leq C_{4}(\Omega_{a_{k}}^2\varepsilon_{a_{k}}^6)^{-\frac{1}{5}}
    \ \ \hbox{as}\ \ k\to\infty,
\end{aligned}
\end{equation*}
which and \eqref{6.38} imply that
\begin{equation*}
\begin{aligned}
  |\nu_{a_{k}}(x)|=&\big|\bar{\eps}_{a_{k}}\sqrt{a_{k}}u_{a_{k}}(\bar{\eps}_{a_{k}}x+x_{a_{k}})e^{-i(\bar{\eps}_{a_{k}}\Omega_{a_{k}} x\cdot x_{a_{k}}^{\bot}-\theta_{a_{k}})}\big|>0\\
 \qquad\hbox{in}\, \ &\big\{x\in\R^2:\, |x|\leq C_{4}(\Omega_{a_{k}}^2\eps_{a_{k}}^6)^{-\frac{1}{5}}\big\}\ \ \hbox{as}\ \ k\to\infty.
\end{aligned}
\end{equation*}
 Therefore, there exist constants $C_{5}>0$ and $C_{*}>0$ such that
\begin{equation}\label{11}
    |u_{a_{k}}(y)|>0 \ \ \hbox{if}\ \  |y|\leq C_{5}\eps_{a_{k}}(\Omega_{a_{k}}^2\varepsilon_{a_{k}}^6)^{-\frac{1}{5}}\leq C_{*}(a^*-a_{k})^{-\frac{1-6\beta}{20}}\ \ \hbox{as}\ \ k\to\infty,
\end{equation}
where $\beta\in[0,\frac{1}{6})$, and $\eps_{a_{k}}=\frac{1}{\sqrt{C_{0}}\lambda}(a^*-a_{k})^{\frac{1+2\beta}{4}}>0$.
We then deduce that $u_{a_{k}}$ does not admit any vortex in the region $R_{a_{k}}:=\{x\in\R^2:\ |x|\leq C_{*}(a^*-a_{k})^{-\frac{1-6\beta}{20}}\}$
as $k\to\infty$, where $\beta\in [0,\frac{1}{6})$ and $C_{*}>0$ is small and independent of $k>0$. In fact, \eqref{11} holds for any sequence $\{a_{k}\}$ satisfying $a_{k}\nearrow a^*$ as $k\to\infty$, which completes the proof of Theorem \ref{th4}.
\qed

\end{document}